

\mag 1200
\input amstex

\let\USversion\relax

\expandafter\ifx\csname selberg.def\endcsname\relax \else\endinput\fi
\expandafter\edef\csname selberg.def\endcsname{%
 \catcode`\noexpand\@=\the\catcode`\@\space}

\let\atbefore @

\catcode`\@=11

\overfullrule\z@

\def\PaperA4{\hsize 6.25truein \vsize 9.63truein}

\def\PaperUS{\hsize 6.6truein \vsize 9truein} 

\def\foliorm{\ifMag\eightrm\else\ninerm\fi}

\let\@ft@\expandafter \let\@tb@f@\atbefore

\newif\ifMag
\def\Magset{\ifnum\mag>\@m\Magtrue\fi}
\Magset

\newif\ifUS

\newdimen\p@@ \p@@\p@
\def\m@ths@r{\ifnum\mathsurround=\z@\z@\else\maths@r\fi}
\def\maths@r{1.6\p@@} \def\mathsurzero{\def\maths@r{\z@}}

\mathsurround\maths@r
\font\Brm=cmr12 \font\Bbf=cmbx12 \font\Bit=cmti12 \font\ssf=cmss10
\font\Bsl=cmsl10 scaled 1200 \font\Bmmi=cmmi10 scaled 1200
\font\BBf=cmbx12 scaled 1200 \font\BMmi=cmmi10 scaled 1440

\def\atletter{\edef\atrestore{\catcode`\noexpand\@=\the\catcode`\@\space}
 \catcode`\@=11}

\newread\@ux \newwrite\@@x \newwrite\@@cd
\let\@np@@\input
\def\@np@t#1{\openin\@ux#1\relax\ifeof\@ux\else\closein\@ux\relax\@np@@ #1\fi}
\def\input#1 {\openin\@ux#1\relax\ifeof\@ux\wrs@x{No file #1}\else
 \closein\@ux\relax\@np@@ #1\fi}
\def\Input#1 {\relax} 

\def\wr@@x#1{} \def\wrs@x{\immediate\write\sixt@@n}

\def\readldf{\@np@t{\jobname.ldf}}
\def\writeldf{\def\wr@@x{\immediate\write\@@x}\def\wr@x@{\write\@@x}
 \def\cl@selbl{\wr@@x{\string\def\string\nextpage{\the\pageno}}%
 \wr@@x{\string\endinput}\immediate\closeout\@@x}
 \immediate\openout\@@x\jobname.ldf}
\let\cl@selbl\relax

\def\nextpage{1}

\def\tod@y{\ifcase\month\or
 January\or February\or March\or April\or May\or June\or July\or
 August\or September\or October\or November\or December\fi\space\,
\number\day,\space\,\number\year}

\newcount\c@time
\def\h@@r{hh}\def\m@n@te{mm}
\def\wh@tt@me{\c@time\time\divide\c@time 60\edef\h@@r{\number\c@time}%
 \multiply\c@time -60\advance\c@time\time\edef
 \m@n@te{\ifnum\c@time<10 0\fi\number\c@time}}
\def\t@me{\h@@r\/{\rm:}\m@n@te} \let\whattime\wh@tt@me
\def\today{\tod@y\wr@@x{\string\todaydef{\tod@y}}}
\def\nowtime{\t@me{\let\/\ic@\wr@@x{\string\nowtimedef{\t@me}}}}
\def\todaydef#1{} \def\nowtimedef#1{}

 \def\emph#1{{\sl #1\/}}

\def\fitem#1{\par\setbox\z@\hbox{#1}\hangindent\wd\z@
 \hglue-2\parindent\kern\wd\z@\indent\llap{#1}\ignore}

\def\itemflat#1{\par\setbox\z@\hbox{\rm #1\enspace}\hang\ifnum\wd\z@>\parindent
 \noindent\unhbox\z@\ignore\else\textindent{\rm#1}\fi}

\newcount\itemlet
\def\newbi{\itemlet 96} \newbi
\def\bitem{\gad\itemlet \par\hangindent1.5\parindent
 \hglue-.5\parindent\textindent{\rm\rlap{\char\the\itemlet}\hp{b})}}

\newcount\itemrm

\def\iitem{\gad\itemrm \par\hangindent1.5\parindent
 \hglue-.5\parindent\textindent{\rm\hp{v}\llap{\romannumeral\the\itemrm})}}

\newcount\itemar

\def\iitema{\gad\itemrm \par\hangindent1.5\parindent
 \hglue-.5\parindent\textindent{\rm\hp{0}\llap{\the\itemrm}.}}

\def\center{\par\begingroup\leftskip\z@ plus \hsize \rightskip\leftskip
 \parindent\z@\parfillskip\z@skip \def\\{\unskip\break}}
\def\endcenter{\endgraf\endgroup}

\def\Abstract{\begingroup\narrower\nt{\bf Abstract.}\enspace\ignore}
\def\endAbs{\endgraf\endgroup}

\let\b@gr@@\begingroup \let\B@gr@@\begingroup
\def\b@gr@{\b@gr@@\let\b@gr@@\undefined}
\def\B@gr@{\B@gr@@\let\B@gr@@\undefined}

\def\@fn@xt#1#2#3{\let\@ch@r=#1\def\n@xt{\ifx\t@st@\@ch@r
 \def\n@@xt{#2}\else\def\n@@xt{#3}\fi\n@@xt}\futurelet\t@st@\n@xt}

\def\@fwd@@#1#2#3{\setbox\z@\hbox{#1}\ifdim\wd\z@>\z@#2\else#3\fi}
\def\s@twd@#1#2{\setbox\z@\hbox{#2}#1\wd\z@}

\def\r@st@re#1{\let#1\s@v@} \def\s@v@d@f{\let\s@v@}

\def\p@sk@p#1#2{\par\skip@#2\relax\ifdim\lastskip<\skip@\relax\removelastskip
 \ifnum#1=\z@\else\penalty#1\relax\fi\vskip\skip@
 \else\ifnum#1=\z@\else\penalty#1\relax\fi\fi}
\def\sk@@p#1{\par\skip@#1\relax\ifdim\lastskip<\skip@\relax\removelastskip
 \vskip\skip@\fi}

\newbox\p@b@ld
\def\poorbold#1{\setbox\p@b@ld\hbox{#1}\kern-.01em\copy\p@b@ld\kern-\wd\p@b@ld
 \kern.02em\copy\p@b@ld\kern-\wd\p@b@ld\kern-.012em\raise.02em\box\p@b@ld}

\ifx\plainfootnote\undefined \let\plainfootnote\footnote \fi

\let\s@v@\proclaim \let\proclaim\relax
\def\r@R@fs#1{\let#1\s@R@fs} \let\s@R@fs\Refs \let\Refs\relax
\def\r@endd@#1{\let#1\s@endd@} \let\s@endd@\enddocument
\let\bye\relax

\def\myR@fs{\@fn@xt[\m@R@f@\m@R@fs} \def\m@R@fs{\@fn@xt*\m@r@f@@\m@R@f@@}
\def\m@R@f@@{\m@R@f@[References]} \def\m@r@f@@*{\m@R@f@[]}

\def\Twelvepoint{\twelvepoint \let\Bbf\BBf \let\Bmmi\BMmi
\font\Brm=cmr12 scaled 1200 \font\Bit=cmti12 scaled 1200
\font\ssf=cmss10 scaled 1200 \font\Bsl=cmsl10 scaled 1440
\font\BBf=cmbx12 scaled 1440 \font\BMmi=cmmi10 scaled 1728}

\newdimen\b@gsize

\newdimen\r@f@nd \newbox\r@f@b@x \newbox\adjb@x
\newbox\p@nct@ \newbox\k@yb@x \newcount\rcount
\newbox\b@b@x \newbox\p@p@rb@x \newbox\j@@rb@x \newbox\y@@rb@x
\newbox\v@lb@x \newbox\is@b@x \newbox\p@g@b@x \newif\ifp@g@ \newif\ifp@g@s
\newbox\inb@@kb@x \newbox\b@@kb@x \newbox\p@blb@x \newbox\p@bl@db@x
\newbox\ed@b@x \newif\ifed@ \newif\ifed@s \newif\if@fl@b \newif\if@fn@m
\newbox\p@p@nf@b@x \newbox\inf@b@x \newbox\b@@nf@b@x
\newtoks\@dd@p@n \newtoks\@ddt@ks

\newif\ifp@gen@

\def\p@@nt{.\kern.3em} \let\point\p@@nt

\let\proheadfont\bf \let\probodyfont\sl \let\demofont\it

\headline={\hfil}
\footline={\ifp@gen@\ifnum\pageno=\z@\else\hfil\foliorm\folio\fi\else
 \ifnum\pageno=\z@\hfil\foliorm\folio\fi\fi\hfil\global\p@gen@true}
\parindent1pc

\font@\tensmc=cmcsc10
\font@\sevenex=cmex7
\font@\sevenit=cmti7
\font@\eightrm=cmr8
\font@\sixrm=cmr6
\font@\eighti=cmmi8 \skewchar\eighti='177
\font@\sixi=cmmi6 \skewchar\sixi='177
\font@\eightsy=cmsy8 \skewchar\eightsy='60
\font@\sixsy=cmsy6 \skewchar\sixsy='60
\font@\eightex=cmex8
\font@\eightbf=cmbx8
\font@\sixbf=cmbx6
\font@\eightit=cmti8
\font@\eightsl=cmsl8
\font@\eightsmc=cmcsc8
\font@\eighttt=cmtt8
\font@\ninerm=cmr9
\font@\ninei=cmmi9 \skewchar\ninei='177
\font@\ninesy=cmsy9 \skewchar\ninesy='60
\font@\nineex=cmex9
\font@\ninebf=cmbx9
\font@\nineit=cmti9
\font@\ninesl=cmsl9
\font@\ninesmc=cmcsc9
\font@\ninemsa=msam9
\font@\ninemsb=msbm9
\font@\nineeufm=eufm9
\font@\eightmsa=msam8
\font@\eightmsb=msbm8
\font@\eighteufm=eufm8
\font@\sixmsa=msam6
\font@\sixmsb=msbm6
\font@\sixeufm=eufm6

\loadmsam\loadmsbm\loadeufm
\input amssym.tex

\def\footnoterule{\kern-3\p@\hrule width5pc\kern 2.6\p@}
\def\m@k@foot#1{\insert\footins
 {\interlinepenalty\interfootnotelinepenalty
 \ifMag\eightpoint\else\ninepoint\fi
 \splittopskip\ht\strutbox\splitmaxdepth\dp\strutbox
 \floatingpenalty\@MM\leftskip\z@\rightskip\z@
 \spaceskip\z@\xspaceskip\z@
 \leavevmode\footstrut\ignore#1\unskip\lower\dp\strutbox
 \vbox to\dp\strutbox{}}}
\def\ftext#1{\m@k@foot{\vsk-.8>\nt #1}}
\def\pr@cl@@m#1{\p@sk@p{-100}\medskipamount
 \def\endproclaim{\endgroup\p@sk@p{55}\medskipamount}\begingroup
 \nt\ignore\proheadfont#1\unskip.\enspace\probodyfont\ignore}
\outer\def\proclaim{\pr@cl@@m} \s@v@d@f\proclaim \let\proclaim\relax
\def\demo#1{\sk@@p\medskipamount\nt{\ignore\demofont#1\unskip.}\enspace
 \ignore}
\def\enddemo{\sk@@p\medskipamount}

\def\cite#1{{\rm[#1]}} 
 \def\Refs#1#2{\relax}

\def\big@#1#2{{\hbox{$\left#2\vcenter to#1\b@gsize{}%
 \right.\nulldelimiterspace\z@\m@th$}}}
\def\big{\big@\@ne}
\def\Big{\big@{1.5}}
\def\bigg{\big@\tw@}
\def\Bigg{\big@{2.5}}
\normallineskiplimit\p@

\def\tenpoint{\p@@\p@ \normallineskiplimit\p@@
 \mathsurround\m@ths@r \normalbaselineskip12\p@@
 \abovedisplayskip12\p@@ plus3\p@@ minus9\p@@
 \belowdisplayskip\abovedisplayskip
 \abovedisplayshortskip\z@ plus3\p@@
 \belowdisplayshortskip7\p@@ plus3\p@@ minus4\p@@
 \textonlyfont@\rm\tenrm \textonlyfont@\it\tenit
 \textonlyfont@\sl\tensl \textonlyfont@\bf\tenbf
 \textonlyfont@\smc\tensmc \textonlyfont@\tt\tentt
 \ifsyntax@ \def\big##1{{\hbox{$\left##1\right.$}}}%
  \let\Big\big \let\bigg\big \let\Bigg\big
 \else
  \textfont\z@\tenrm \scriptfont\z@\sevenrm \scriptscriptfont\z@\fiverm
  \textfont\@ne\teni \scriptfont\@ne\seveni \scriptscriptfont\@ne\fivei
  \textfont\tw@\tensy \scriptfont\tw@\sevensy \scriptscriptfont\tw@\fivesy
  \textfont\thr@@\tenex \scriptfont\thr@@\sevenex
	\scriptscriptfont\thr@@\sevenex
  \textfont\itfam\tenit \scriptfont\itfam\sevenit
	\scriptscriptfont\itfam\sevenit
  \textfont\bffam\tenbf \scriptfont\bffam\sevenbf
	\scriptscriptfont\bffam\fivebf
  \textfont\msafam\tenmsa \scriptfont\msafam\sevenmsa
	\scriptscriptfont\msafam\fivemsa
  \textfont\msbfam\tenmsb \scriptfont\msbfam\sevenmsb
	\scriptscriptfont\msbfam\fivemsb
  \textfont\eufmfam\teneufm \scriptfont\eufmfam\seveneufm
	\scriptscriptfont\eufmfam\fiveeufm
  \setbox\strutbox\hbox{\vrule height8.5\p@@ depth3.5\p@@ width\z@}%
  \setbox\strutbox@\hbox{\lower.5\normallineskiplimit\vbox{%
	\kern-\normallineskiplimit\copy\strutbox}}%
   \setbox\z@\vbox{\hbox{$($}\kern\z@}\b@gsize1.2\ht\z@
  \fi
  \normalbaselines\rm\dotsspace@1.5mu\ex@.2326ex\jot3\ex@}

\def\eightpoint{\p@@.8\p@ \normallineskiplimit\p@@
 \mathsurround\m@ths@r \normalbaselineskip10\p@
 \abovedisplayskip10\p@ plus2.4\p@ minus7.2\p@
 \belowdisplayskip\abovedisplayskip
 \abovedisplayshortskip\z@ plus3\p@@
 \belowdisplayshortskip7\p@@ plus3\p@@ minus4\p@@
 \textonlyfont@\rm\eightrm \textonlyfont@\it\eightit
 \textonlyfont@\sl\eightsl \textonlyfont@\bf\eightbf
 \textonlyfont@\smc\eightsmc \textonlyfont@\tt\eighttt
 \ifsyntax@\def\big##1{{\hbox{$\left##1\right.$}}}%
  \let\Big\big \let\bigg\big \let\Bigg\big
 \else
  \textfont\z@\eightrm \scriptfont\z@\sixrm \scriptscriptfont\z@\fiverm
  \textfont\@ne\eighti \scriptfont\@ne\sixi \scriptscriptfont\@ne\fivei
  \textfont\tw@\eightsy \scriptfont\tw@\sixsy \scriptscriptfont\tw@\fivesy
  \textfont\thr@@\eightex \scriptfont\thr@@\sevenex
	\scriptscriptfont\thr@@\sevenex
  \textfont\itfam\eightit \scriptfont\itfam\sevenit
	\scriptscriptfont\itfam\sevenit
  \textfont\bffam\eightbf \scriptfont\bffam\sixbf
	\scriptscriptfont\bffam\fivebf
  \textfont\msafam\eightmsa \scriptfont\msafam\sixmsa
	\scriptscriptfont\msafam\fivemsa
  \textfont\msbfam\eightmsb \scriptfont\msbfam\sixmsb
	\scriptscriptfont\msbfam\fivemsb
  \textfont\eufmfam\eighteufm \scriptfont\eufmfam\sixeufm
	\scriptscriptfont\eufmfam\fiveeufm
 \setbox\strutbox\hbox{\vrule height7\p@ depth3\p@ width\z@}%
 \setbox\strutbox@\hbox{\raise.5\normallineskiplimit\vbox{%
   \kern-\normallineskiplimit\copy\strutbox}}%
 \setbox\z@\vbox{\hbox{$($}\kern\z@}\b@gsize1.2\ht\z@
 \fi
 \normalbaselines\eightrm\dotsspace@1.5mu\ex@.2326ex\jot3\ex@}

\def\ninepoint{\p@@.9\p@ \normallineskiplimit\p@@
 \mathsurround\m@ths@r \normalbaselineskip11\p@
 \abovedisplayskip11\p@ plus2.7\p@ minus8.1\p@
 \belowdisplayskip\abovedisplayskip
 \abovedisplayshortskip\z@ plus3\p@@
 \belowdisplayshortskip7\p@@ plus3\p@@ minus4\p@@
 \textonlyfont@\rm\ninerm \textonlyfont@\it\nineit
 \textonlyfont@\sl\ninesl \textonlyfont@\bf\ninebf
 \textonlyfont@\smc\ninesmc \textonlyfont@\tt\ninett
 \ifsyntax@ \def\big##1{{\hbox{$\left##1\right.$}}}%
  \let\Big\big \let\bigg\big \let\Bigg\big
 \else
  \textfont\z@\ninerm \scriptfont\z@\sevenrm \scriptscriptfont\z@\fiverm
  \textfont\@ne\ninei \scriptfont\@ne\seveni \scriptscriptfont\@ne\fivei
  \textfont\tw@\ninesy \scriptfont\tw@\sevensy \scriptscriptfont\tw@\fivesy
  \textfont\thr@@\nineex \scriptfont\thr@@\sevenex
	\scriptscriptfont\thr@@\sevenex
  \textfont\itfam\nineit \scriptfont\itfam\sevenit
	\scriptscriptfont\itfam\sevenit
  \textfont\bffam\ninebf \scriptfont\bffam\sevenbf
	\scriptscriptfont\bffam\fivebf
  \textfont\msafam\ninemsa \scriptfont\msafam\sevenmsa
	\scriptscriptfont\msafam\fivemsa
  \textfont\msbfam\ninemsb \scriptfont\msbfam\sevenmsb
	\scriptscriptfont\msbfam\fivemsb
  \textfont\eufmfam\nineeufm \scriptfont\eufmfam\seveneufm
	\scriptscriptfont\eufmfam\fiveeufm
  \setbox\strutbox\hbox{\vrule height8.5\p@@ depth3.5\p@@ width\z@}%
  \setbox\strutbox@\hbox{\lower.5\normallineskiplimit\vbox{%
	\kern-\normallineskiplimit\copy\strutbox}}%
   \setbox\z@\vbox{\hbox{$($}\kern\z@}\b@gsize1.2\ht\z@
  \fi
  \normalbaselines\rm\dotsspace@1.5mu\ex@.2326ex\jot3\ex@}

\font@\twelverm=cmr10 scaled 1200
\font@\twelveit=cmti10 scaled 1200
\font@\twelvesl=cmsl10 scaled 1200
\font@\twelvebf=cmbx10 scaled 1200
\font@\twelvesmc=cmcsc10 scaled 1200
\font@\twelvett=cmtt10 scaled 1200
\font@\twelvei=cmmi10 scaled 1200 \skewchar\twelvei='177
\font@\twelvesy=cmsy10 scaled 1200 \skewchar\twelvesy='60
\font@\twelveex=cmex10 scaled 1200
\font@\twelvemsa=msam10 scaled 1200
\font@\twelvemsb=msbm10 scaled 1200
\font@\twelveeufm=eufm10 scaled 1200

\def\twelvepoint{\p@@1.2\p@ \normallineskiplimit\p@@
 \mathsurround\m@ths@r \normalbaselineskip12\p@@
 \abovedisplayskip12\p@@ plus3\p@@ minus9\p@@
 \belowdisplayskip\abovedisplayskip
 \abovedisplayshortskip\z@ plus3\p@@
 \belowdisplayshortskip7\p@@ plus3\p@@ minus4\p@@
 \textonlyfont@\rm\twelverm \textonlyfont@\it\twelveit
 \textonlyfont@\sl\twelvesl \textonlyfont@\bf\twelvebf
 \textonlyfont@\smc\twelvesmc \textonlyfont@\tt\twelvett
 \ifsyntax@ \def\big##1{{\hbox{$\left##1\right.$}}}%
  \let\Big\big \let\bigg\big \let\Bigg\big
 \else
  \textfont\z@\twelverm \scriptfont\z@\eightrm \scriptscriptfont\z@\sixrm
  \textfont\@ne\twelvei \scriptfont\@ne\eighti \scriptscriptfont\@ne\sixi
  \textfont\tw@\twelvesy \scriptfont\tw@\eightsy \scriptscriptfont\tw@\sixsy
  \textfont\thr@@\twelveex \scriptfont\thr@@\eightex
	\scriptscriptfont\thr@@\sevenex
  \textfont\itfam\twelveit \scriptfont\itfam\eightit
	\scriptscriptfont\itfam\sevenit
  \textfont\bffam\twelvebf \scriptfont\bffam\eightbf
	\scriptscriptfont\bffam\sixbf
  \textfont\msafam\twelvemsa \scriptfont\msafam\eightmsa
	\scriptscriptfont\msafam\sixmsa
  \textfont\msbfam\twelvemsb \scriptfont\msbfam\eightmsb
	\scriptscriptfont\msbfam\sixmsb
  \textfont\eufmfam\twelveeufm \scriptfont\eufmfam\eighteufm
	\scriptscriptfont\eufmfam\sixeufm
  \setbox\strutbox\hbox{\vrule height8.5\p@@ depth3.5\p@@ width\z@}%
  \setbox\strutbox@\hbox{\lower.5\normallineskiplimit\vbox{%
	\kern-\normallineskiplimit\copy\strutbox}}%
  \setbox\z@\vbox{\hbox{$($}\kern\z@}\b@gsize1.2\ht\z@
  \fi
  \normalbaselines\rm\dotsspace@1.5mu\ex@.2326ex\jot3\ex@}

\font@\twelvetrm=cmr10 at 12truept
\font@\twelvetit=cmti10 at 12truept
\font@\twelvetsl=cmsl10 at 12truept
\font@\twelvetbf=cmbx10 at 12truept
\font@\twelvetsmc=cmcsc10 at 12truept
\font@\twelvettt=cmtt10 at 12truept
\font@\twelveti=cmmi10 at 12truept \skewchar\twelveti='177
\font@\twelvetsy=cmsy10 at 12truept \skewchar\twelvetsy='60
\font@\twelvetex=cmex10 at 12truept
\font@\twelvetmsa=msam10 at 12truept
\font@\twelvetmsb=msbm10 at 12truept
\font@\twelveteufm=eufm10 at 12truept

\def\twelvetruepoint{\p@@1.2truept \normallineskiplimit\p@@
 \mathsurround\m@ths@r \normalbaselineskip12\p@@
 \abovedisplayskip12\p@@ plus3\p@@ minus9\p@@
 \belowdisplayskip\abovedisplayskip
 \abovedisplayshortskip\z@ plus3\p@@
 \belowdisplayshortskip7\p@@ plus3\p@@ minus4\p@@
 \textonlyfont@\rm\twelvetrm \textonlyfont@\it\twelvetit
 \textonlyfont@\sl\twelvetsl \textonlyfont@\bf\twelvetbf
 \textonlyfont@\smc\twelvetsmc \textonlyfont@\tt\twelvettt
 \ifsyntax@ \def\big##1{{\hbox{$\left##1\right.$}}}%
  \let\Big\big \let\bigg\big \let\Bigg\big
 \else
  \textfont\z@\twelvetrm \scriptfont\z@\eightrm \scriptscriptfont\z@\sixrm
  \textfont\@ne\twelveti \scriptfont\@ne\eighti \scriptscriptfont\@ne\sixi
  \textfont\tw@\twelvetsy \scriptfont\tw@\eightsy \scriptscriptfont\tw@\sixsy
  \textfont\thr@@\twelvetex \scriptfont\thr@@\eightex
	\scriptscriptfont\thr@@\sevenex
  \textfont\itfam\twelvetit \scriptfont\itfam\eightit
	\scriptscriptfont\itfam\sevenit
  \textfont\bffam\twelvetbf \scriptfont\bffam\eightbf
	\scriptscriptfont\bffam\sixbf
  \textfont\msafam\twelvetmsa \scriptfont\msafam\eightmsa
	\scriptscriptfont\msafam\sixmsa
  \textfont\msbfam\twelvetmsb \scriptfont\msbfam\eightmsb
	\scriptscriptfont\msbfam\sixmsb
  \textfont\eufmfam\twelveteufm \scriptfont\eufmfam\eighteufm
	\scriptscriptfont\eufmfam\sixeufm
  \setbox\strutbox\hbox{\vrule height8.5\p@@ depth3.5\p@@ width\z@}%
  \setbox\strutbox@\hbox{\lower.5\normallineskiplimit\vbox{%
	\kern-\normallineskiplimit\copy\strutbox}}%
  \setbox\z@\vbox{\hbox{$($}\kern\z@}\b@gsize1.2\ht\z@
  \fi
  \normalbaselines\rm\dotsspace@1.5mu\ex@.2326ex\jot3\ex@}

\font@\elevenrm=cmr10 scaled 1095
\font@\elevenit=cmti10 scaled 1095
\font@\elevensl=cmsl10 scaled 1095
\font@\elevenbf=cmbx10 scaled 1095
\font@\elevensmc=cmcsc10 scaled 1095
\font@\eleventt=cmtt10 scaled 1095
\font@\eleveni=cmmi10 scaled 1095 \skewchar\eleveni='177
\font@\elevensy=cmsy10 scaled 1095 \skewchar\elevensy='60
\font@\elevenex=cmex10 scaled 1095
\font@\elevenmsa=msam10 scaled 1095
\font@\elevenmsb=msbm10 scaled 1095
\font@\eleveneufm=eufm10 scaled 1095

\def\elevenpoint{\p@@1.1\p@ \normallineskiplimit\p@@
 \mathsurround\m@ths@r \normalbaselineskip12\p@@
 \abovedisplayskip12\p@@ plus3\p@@ minus9\p@@
 \belowdisplayskip\abovedisplayskip
 \abovedisplayshortskip\z@ plus3\p@@
 \belowdisplayshortskip7\p@@ plus3\p@@ minus4\p@@
 \textonlyfont@\rm\elevenrm \textonlyfont@\it\elevenit
 \textonlyfont@\sl\elevensl \textonlyfont@\bf\elevenbf
 \textonlyfont@\smc\elevensmc \textonlyfont@\tt\eleventt
 \ifsyntax@ \def\big##1{{\hbox{$\left##1\right.$}}}%
  \let\Big\big \let\bigg\big \let\Bigg\big
 \else
  \textfont\z@\elevenrm \scriptfont\z@\eightrm \scriptscriptfont\z@\sixrm
  \textfont\@ne\eleveni \scriptfont\@ne\eighti \scriptscriptfont\@ne\sixi
  \textfont\tw@\elevensy \scriptfont\tw@\eightsy \scriptscriptfont\tw@\sixsy
  \textfont\thr@@\elevenex \scriptfont\thr@@\eightex
	\scriptscriptfont\thr@@\sevenex
  \textfont\itfam\elevenit \scriptfont\itfam\eightit
	\scriptscriptfont\itfam\sevenit
  \textfont\bffam\elevenbf \scriptfont\bffam\eightbf
	\scriptscriptfont\bffam\sixbf
  \textfont\msafam\elevenmsa \scriptfont\msafam\eightmsa
	\scriptscriptfont\msafam\sixmsa
  \textfont\msbfam\elevenmsb \scriptfont\msbfam\eightmsb
	\scriptscriptfont\msbfam\sixmsb
  \textfont\eufmfam\eleveneufm \scriptfont\eufmfam\eighteufm
	\scriptscriptfont\eufmfam\sixeufm
  \setbox\strutbox\hbox{\vrule height8.5\p@@ depth3.5\p@@ width\z@}%
  \setbox\strutbox@\hbox{\lower.5\normallineskiplimit\vbox{%
	\kern-\normallineskiplimit\copy\strutbox}}%
  \setbox\z@\vbox{\hbox{$($}\kern\z@}\b@gsize1.2\ht\z@
  \fi
  \normalbaselines\rm\dotsspace@1.5mu\ex@.2326ex\jot3\ex@}

\def\m@R@f@[#1]{\mathsurzero{
 \s@ct{}{#1}}\wr@@c{\string\Refcd{#1}{\the\pageno}}\B@gr@
 \frenchspacing\rcount\z@\refkey{\k@yf@nt[##1]}\refno{\k@yf@nt[##1]}%
 \widest{AZ}\keyright\let\Key\key\let\refin\relax}
\def\widest#1{\s@twd@\r@f@nd{\r@fk@y{\k@yf@nt#1}\enspace}}
\def\widestno#1{\s@twd@\r@f@nd{\r@fn@{\k@yf@nt#1}\enspace}}
\def\widestlabel#1{\s@twd@\r@f@nd{\k@yf@nt#1\enspace}}
\def\refkey{\def\r@fk@y##1} \def\refno{\def\r@fn@##1}
\def\keyright{\def\r@fit@m{\hang\textindent}}
\def\keyflat{\def\r@fit@m##1{\setbox\z@\hbox{##1\enspace}\hang\noindent
 \ifnum\wd\z@<\parindent\indent\hglue-\wd\z@\fi\unhbox\z@}}

\def\R@fb@x{\global\setbox\r@f@b@x} \def\K@yb@x{\global\setbox\k@yb@x}
\def\ref{\par\b@gr@\r@ff@nt\R@fb@x\box\voidb@x\K@yb@x\box\voidb@x
 \@fn@mfalse\@fl@bfalse\b@g@nr@f}
\def\c@nc@t#1{\setbox\z@\lastbox
 \setbox\adjb@x\hbox{\unhbox\adjb@x\unhbox\z@\unskip\unskip\unpenalty#1}}
\def\adjust#1{\relax\ifmmode\penalty-\@M\null\hfil$\clubpenalty\z@
 \widowpenalty\z@\interlinepenalty\z@\offinterlineskip\endgraf
 \setbox\z@\lastbox\unskip\unpenalty\c@nc@t{#1}\nt$\hfil\penalty-\@M
 \else\endgraf\c@nc@t{#1}\nt\fi}
\def\adjustnext#1{\P@nct\hbox{#1}\ignore}
\def\adjustend#1{\def\@djp@{#1}\ignore}
\def\addtoks#1{\global\@ddt@ks{#1}\ignore}
\def\addnext#1{\global\@dd@p@n{#1}\ignore}

\def\cl@s@{\adjust{\@djp@}\endgraf\setbox\z@\lastbox
 \global\setbox\@ne\hbox{\unhbox\adjb@x\ifvoid\z@\else\unhbox\z@\unskip\unskip
 \unpenalty\fi}\egroup\ifnum\c@rr@nt=\k@yb@x\global\fi
 \setbox\c@rr@nt\hbox{\unhbox\@ne\box\p@nct@}\P@nct\null
 \the\@ddt@ks\global\@ddt@ks{}}
\def\@p@n#1{\def\c@rr@nt{#1}\setbox\c@rr@nt\vbox\bgroup\let\@djp@\relax
 \hsize\maxdimen\nt\the\@dd@p@n\global\@dd@p@n{}}
\def\b@g@nr@f{\bgroup\@p@n\z@}
\def\key{\cl@s@\ifvoid\k@yb@x\@p@n\k@yb@x\k@yf@nt\else\@p@n\z@\fi}
\def\label{\cl@s@\ifvoid\k@yb@x\global\@fl@btrue\@p@n\k@yb@x\k@yf@nt\else
 \@p@n\z@\fi}
\def\no{\cl@s@\ifvoid\k@yb@x\gad\rcount\global\@fn@mtrue
 \K@yb@x\hbox{\k@yf@nt\the\rcount}\fi\@p@n\z@}
\def\labelno{\cl@s@\ifvoid\k@yb@x\gad\rcount\@fl@btrue
 \@p@n\k@yb@x\k@yf@nt\the\rcount\else\@p@n\z@\fi}
\def\by{\cl@s@\@p@n\b@b@x} \def\paper{\cl@s@\@p@n\p@p@rb@x\p@p@rf@nt\ignore}
\def\jour{\cl@s@\@p@n\j@@rb@x} \def\yr{\cl@s@\@p@n\y@@rb@x}
\def\vol{\cl@s@\@p@n\v@lb@x\v@lf@nt\ignore}
\def\issue{\cl@s@\@p@n\is@b@x\iss@f@nt\ignore}
\def\page{\cl@s@\ifp@g@s\@p@n\z@\else\p@g@true\@p@n\p@g@b@x\fi}
\def\pages{\cl@s@\ifp@g@\@p@n\z@\else\p@g@strue\@p@n\p@g@b@x\fi}
\def\inbook{\cl@s@\@p@n\inb@@kb@x}
\def\book{\cl@s@\@p@n\b@@kb@x\b@@kf@nt\ignore}
\def\publ{\cl@s@\@p@n\p@blb@x} \def\publaddr{\cl@s@\@p@n\p@bl@db@x}
\def\ed{\cl@s@\ifed@s\@p@n\z@\else\ed@true\@p@n\ed@b@x\fi}
\def\eds{\cl@s@\ifed@\@p@n\z@\else\ed@strue\@p@n\ed@b@x\fi}
\def\info{\cl@s@\@p@n\inf@b@x} \def\paperinfo{\cl@s@\@p@n\p@p@nf@b@x}
\def\bookinfo{\cl@s@\@p@n\b@@nf@b@x} 
\def\P@nct{\global\setbox\p@nct@} \def\nopunct{\P@nct\box\voidb@x}
\def\p@@@t#1#2{\ifvoid\p@nct@\else#1\unhbox\p@nct@#2\fi}
\def\sp@@{\penalty-50 \space\hskip\z@ plus.1em}
\def\c@mm@{\p@@@t,\sp@@} \def\sp@c@{\p@@@t\empty\sp@@}
\def\p@tb@x#1#2{\ifvoid#1\else#2\@nb@x#1\fi}
\def\@nb@x#1{\unhbox#1\P@nct\lastbox}
\def\endr@f@{\cl@s@\nopunct
 \R@fb@x\hbox{\unhbox\r@f@b@x \p@tb@x\b@b@x\empty
 \ifvoid\j@@rb@x\ifvoid\inb@@kb@x\ifvoid\p@p@rb@x\ifvoid\b@@kb@x
  \ifvoid\p@p@nf@b@x\ifvoid\b@@nf@b@x
  \p@tb@x\v@lb@x\c@mm@ \ifvoid\y@@rb@x\else\sp@c@(\@nb@x\y@@rb@x)\fi
  \p@tb@x\is@b@x\c@mm@ \p@tb@x\p@g@b@x\c@mm@ \p@tb@x\inf@b@x\c@mm@
  \else\p@tb@x \b@@nf@b@x\c@mm@ \p@tb@x\v@lb@x\c@mm@ \p@tb@x\is@b@x\sp@c@
  \ifvoid\ed@b@x\else\sp@c@(\@nb@x\ed@b@x,\space\ifed@ ed.\else eds.\fi)\fi
  \p@tb@x\p@blb@x\c@mm@ \p@tb@x\p@bl@db@x\c@mm@ \p@tb@x\y@@rb@x\c@mm@
  \p@tb@x\p@g@b@x{\c@mm@\ifp@g@ p\p@@nt\else pp\p@@nt\fi}%
  \p@tb@x\inf@b@x\c@mm@\fi
  \else \p@tb@x\p@p@nf@b@x\c@mm@ \p@tb@x\v@lb@x\c@mm@
  \ifvoid\y@@rb@x\else\sp@c@(\@nb@x\y@@rb@x)\fi
  \p@tb@x\is@b@x\c@mm@ \p@tb@x\p@g@b@x\c@mm@ \p@tb@x\inf@b@x\c@mm@\fi
  \else \p@tb@x\b@@kb@x\c@mm@
  \p@tb@x\b@@nf@b@x\c@mm@ \p@tb@x\p@blb@x\c@mm@
  \p@tb@x\p@bl@db@x\c@mm@ \p@tb@x\y@@rb@x\c@mm@
  \ifvoid\p@g@b@x\else\c@mm@\@nb@x\p@g@b@x p\fi \p@tb@x\inf@b@x\c@mm@ \fi
  \else \c@mm@\@nb@x\p@p@rb@x\ic@\p@tb@x\p@p@nf@b@x\c@mm@
  \p@tb@x\v@lb@x\sp@c@ \ifvoid\y@@rb@x\else\sp@c@(\@nb@x\y@@rb@x)\fi
  \p@tb@x\is@b@x\c@mm@ \p@tb@x\p@g@b@x\c@mm@\p@tb@x\inf@b@x\c@mm@\fi
  \else \p@tb@x\p@p@rb@x\c@mm@\ic@\p@tb@x\p@p@nf@b@x\c@mm@
  \c@mm@\@nb@x\inb@@kb@x \p@tb@x\b@@nf@b@x\c@mm@ \p@tb@x\v@lb@x\sp@c@
  \p@tb@x\is@b@x\sp@c@
  \ifvoid\ed@b@x\else\sp@c@(\@nb@x\ed@b@x,\space\ifed@ ed.\else eds.\fi)\fi
  \p@tb@x\p@blb@x\c@mm@ \p@tb@x\p@bl@db@x\c@mm@ \p@tb@x\y@@rb@x\c@mm@
  \p@tb@x\p@g@b@x{\c@mm@\ifp@g@ p\p@@nt\else pp\p@@nt\fi}%
  \p@tb@x\inf@b@x\c@mm@\fi
  \else\p@tb@x\p@p@rb@x\c@mm@\ic@\p@tb@x\p@p@nf@b@x\c@mm@\p@tb@x\j@@rb@x\c@mm@
  \p@tb@x\v@lb@x\sp@c@ \ifvoid\y@@rb@x\else\sp@c@(\@nb@x\y@@rb@x)\fi
  \p@tb@x\is@b@x\c@mm@ \p@tb@x\p@g@b@x\c@mm@ \p@tb@x\inf@b@x\c@mm@ \fi}}
\def\m@r@f#1#2{\endr@f@\ifvoid\p@nct@\else\R@fb@x\hbox{\unhbox\r@f@b@x
 #1\unhbox\p@nct@\penalty-200\enskip#2}\fi\egroup\b@g@nr@f}
\def\endref{\endr@f@\ifvoid\p@nct@\else\R@fb@x\hbox{\unhbox\r@f@b@x.}\fi
 \parindent\r@f@nd
 \r@fit@m{\ifvoid\k@yb@x\else\if@fn@m\r@fn@{\unhbox\k@yb@x}\else
 \if@fl@b\unhbox\k@yb@x\else\r@fk@y{\unhbox\k@yb@x}\fi\fi\fi}\unhbox\r@f@b@x
 \endgraf\egroup\endgroup}
\def\moreref{\m@r@f;\empty}
\def\transl{\m@r@f;{\unskip\space
 {\sl English translation\ic@}:\penalty-66 \space}}
\def\endRefs{\endgraf\goodbreak\endgroup}

\hyphenation{acad-e-my acad-e-mies af-ter-thought anom-aly anom-alies
an-ti-deriv-a-tive an-tin-o-my an-tin-o-mies apoth-e-o-ses
apoth-e-o-sis ap-pen-dix ar-che-typ-al as-sign-a-ble as-sist-ant-ship
as-ymp-tot-ic asyn-chro-nous at-trib-uted at-trib-ut-able bank-rupt
bank-rupt-cy bi-dif-fer-en-tial blue-print busier busiest
cat-a-stroph-ic cat-a-stroph-i-cally con-gress cross-hatched data-base
de-fin-i-tive de-riv-a-tive dis-trib-ute dri-ver dri-vers eco-nom-ics
econ-o-mist elit-ist equi-vari-ant ex-quis-ite ex-tra-or-di-nary
flow-chart for-mi-da-ble forth-right friv-o-lous ge-o-des-ic
ge-o-det-ic geo-met-ric griev-ance griev-ous griev-ous-ly
hexa-dec-i-mal ho-lo-no-my ho-mo-thetic ideals idio-syn-crasy
in-fin-ite-ly in-fin-i-tes-i-mal ir-rev-o-ca-ble key-stroke
lam-en-ta-ble light-weight mal-a-prop-ism man-u-script mar-gin-al
meta-bol-ic me-tab-o-lism meta-lan-guage me-trop-o-lis
met-ro-pol-i-tan mi-nut-est mol-e-cule mono-chrome mono-pole
mo-nop-oly mono-spline mo-not-o-nous mul-ti-fac-eted mul-ti-plic-able
non-euclid-ean non-iso-mor-phic non-smooth par-a-digm par-a-bol-ic
pa-rab-o-loid pa-ram-e-trize para-mount pen-ta-gon phe-nom-e-non
post-script pre-am-ble pro-ce-dur-al pro-hib-i-tive pro-hib-i-tive-ly
pseu-do-dif-fer-en-tial pseu-do-fi-nite pseu-do-nym qua-drat-ic
quad-ra-ture qua-si-smooth qua-si-sta-tion-ary qua-si-tri-an-gu-lar
quin-tes-sence quin-tes-sen-tial re-arrange-ment rec-tan-gle
ret-ri-bu-tion retro-fit retro-fit-ted right-eous right-eous-ness
ro-bot ro-bot-ics sched-ul-ing se-mes-ter semi-def-i-nite
semi-ho-mo-thet-ic set-up se-vere-ly side-step sov-er-eign spe-cious
spher-oid spher-oid-al star-tling star-tling-ly sta-tis-tics
sto-chas-tic straight-est strange-ness strat-a-gem strong-hold
sum-ma-ble symp-to-matic syn-chro-nous topo-graph-i-cal tra-vers-a-ble
tra-ver-sal tra-ver-sals treach-ery turn-around un-at-tached
un-err-ing-ly white-space wide-spread wing-spread wretch-ed
wretch-ed-ly Brown-ian Eng-lish Euler-ian Feb-ru-ary Gauss-ian
Grothen-dieck Hamil-ton-ian Her-mit-ian Jan-u-ary Japan-ese Kor-te-weg
Le-gendre Lip-schitz Lip-schitz-ian Mar-kov-ian Noe-ther-ian
No-vem-ber Rie-mann-ian Schwarz-schild Sep-tem-ber}

\let\nopagenumber\p@gen@false \let\putpagenumber\p@gen@true

\outer\def\myRefs{\myR@fs} \r@st@re\proclaim
\def\bye{\par\vfill\supereject\cl@selbl\cl@secd\b@e} \r@endd@\b@e
 \let\Key\key \def\endpro{\par\endproclaim}
\let\d@c@\document \def\document{\d@c@\tenpoint}

\newtoks\@@tp@t \@@tp@t\output
\output=\@ft@{\let\{\noexpand\the\@@tp@t}
\let\{\relax

\newif\ifVersion \Versiontrue
\def\p@n@l#1{\ifnum#1=\z@\else\penalty#1\relax\fi}

\def\s@ct#1#2{\ifVersion
 \skip@\lastskip\ifdim\skip@<1.5\bls\vskip-\skip@\p@n@l{-200}\vsk.5>%
 \p@n@l{-200}\vsk.5>\p@n@l{-200}\vsk.5>\p@n@l{-200}\vsk-1.5>\else
 \p@n@l{-200}\fi\ifdim\skip@<.9\bls\vsk.9>\else
 \ifdim\skip@<1.5\bls\vskip\skip@\fi\fi
 \vtop{\twelvepoint\raggedright\s@cf@nt\vp1\vsk->\vskip.16ex
 \s@twd@\parindent{#1}%
 \ifdim\parindent>\z@\adv\parindent.5em\fi\hang\textindent{#1}#2\strut}
 \else
 \p@sk@p{-200}{.8\bls}\vtop{\s@cf@nt\s@twd@\parindent{#1}%
 \ifdim\parindent>\z@\adv\parindent.5em\fi\hang\textindent{#1}#2\strut}\fi
 \nointerlineskip\nobreak\vtop{\strut}\nobreak\vskip-.6\bls\nobreak}

\def\s@bs@ct#1#2{\ifVersion
 \skip@\lastskip\ifdim\skip@<1.5\bls\vskip-\skip@\p@n@l{-200}\vsk.5>%
 \p@n@l{-200}\vsk.5>\p@n@l{-200}\vsk.5>\p@n@l{-200}\vsk-1.5>\else
 \p@n@l{-200}\fi\ifdim\skip@<.9\bls\vsk.9>\else
 \ifdim\skip@<1.5\bls\vskip\skip@\fi\fi
 \vtop{\elevenpoint\raggedright\s@bf@nt\vp1\vsk->\vskip.16ex%
 \s@twd@\parindent{#1}\ifdim\parindent>\z@\adv\parindent.5em\fi
 \hang\textindent{#1}#2\strut}
 \else
 \p@sk@p{-200}{.6\bls}\vtop{\s@bf@nt\s@twd@\parindent{#1}%
 \ifdim\parindent>\z@\adv\parindent.5em\fi\hang\textindent{#1}#2\strut}\fi
 \nointerlineskip\nobreak\vtop{\strut}\nobreak\vskip-.8\bls\nobreak}

\def\gadv{\global\adv} \def\gad#1{\gadv#1\@ne} \def\gadneg#1{\gadv#1-\@ne}

\newcount\t@@n \t@@n=10 \newbox\testbox

\newcount\Sno \newcount\Lno \newcount\Fno

\def\pr@cl#1{\r@st@re\pr@c@\pr@c@{#1}\global\let\pr@c@\relax}

\def\l@L#1{\l@bel{#1}L} \def\l@F#1{\l@bel{#1}F} \def\<#1>{\l@b@l{#1}F}

\def\tagg#1{\tag"\rlap{\rm(#1)}\kern.01\p@"}
\def\Tag#1{\tag{\l@F{#1}}} \def\Tagg#1{\tagg{\l@F{#1}}}

\def\xspace{\kern.34em}

\def\Th#1{\pr@cl{Theorem\xspace\l@L{#1}}\ignore}
\def\Lm#1{\pr@cl{Lemma\xspace\l@L{#1}}\ignore}
\def\Cr#1{\pr@cl{Corollary\xspace\l@L{#1}}\ignore}
\def\Df#1{\pr@cl{Definition\xspace\l@L{#1}}\ignore}
\def\Cj#1{\pr@cl{Conjecture\xspace\l@L{#1}}\ignore}
\def\Prop#1{\pr@cl{Proposition\xspace\l@L{#1}}\ignore}
 
\def\Pf#1.{\demo{Proof #1}} 

\def\Proof#1.{\demo{\let\{\relax Proof #1}\def\t@st@{#1}%
 \ifx\t@st@\empty\else\xdef\@@wr##1##2##3##4{##1{##2##3{\the\cdn@}{##4}}}%
 \wr@@c@{\the\cdn@}{Proof #1}\@@wr\wr@@c\string\subcd{\the\pageno}\fi\ignore}

\def\Ap@x{Appendix}
\def\Appendix{\Sno=64 \t@@n\@ne \wr@@c{\string\Appencd}
 \def\sf@rm{\char\the\Sno} \def\sf@rm@{\Ap@x\space\sf@rm} \def\sf@rm@@{\Ap@x}
 \def\s@ct@n##1##2{\s@ct\empty{\setbox\z@\hbox{##1}\ifdim\wd\z@=\z@
 \if##2*\sf@rm@@\else\if##2.\sf@rm@@.\else##2\fi\fi\else
 \if##2*\sf@rm@\else\if##2.\sf@rm@.\else\sf@rm@.\enspace##2\fi\fi\fi}}}
\def\Appcd#1#2#3{\gad\Cdentry\global\cdentry\z@\def\Ap@@{\hglue-\l@ftcd\Ap@x}
 \ifx\@ppl@ne\empty\def\l@@b{\@fwd@@{#1}{\space#1}{}}
 \if*#2\entcd{}{\Ap@@\l@@b}{#3}\else\if.#2\entcd{}{\Ap@@\l@@b.}{#3}\else
 \entcd{}{\Ap@@\l@@b.\enspace#2}{#3}\fi\fi\else
 \def\l@@b{\@fwd@@{#1}{\c@l@b{#1}}{}}\if*#2\entcd{\l@@b}{\Ap@x}{#3}\else
 \if.#2\entcd{\l@@b}{\Ap@x.}{#3}\else\entcd{\l@@b}{#2}{#3}\fi\fi\fi}

\let\s@ct@n\s@ct
\def\s@ct@@[#1]#2{\@ft@\xdef\csname @#1@S@\endcsname{\sf@rm}\wr@@x{}%
 \wr@@x{\string\labeldef{S}\space{\?#1@S?}\space{#1}}%
 {
 \s@ct@n{\sf@rm@}{#2}}\wr@@c{\string\Entcd{\?#1@S?}{#2}{\the\pageno}}}
\def\s@ct@#1{\wr@@x{}{
 \s@ct@n{\sf@rm@}{#1}}\wr@@c{\string\Entcd{\sf@rm}{#1}{\the\pageno}}}
\def\s@ct@e[#1]#2{\@ft@\xdef\csname @#1@S@\endcsname{\sf@rm}\wr@@x{}%
 \wr@@x{\string\labeldef{S}\space{\?#1@S?}\space{#1}}%
 {
 \s@ct@n\empty{#2}}\wr@@c{\string\Entcd{}{#2}{\the\pageno}}}
\def\s@cte#1{\wr@@x{}{
 \s@ct@n\empty{#1}}\wr@@c{\string\Entcd{}{#1}{\the\pageno}}}
\def\theSno#1#2{\dff\?#1@S?{#2}%
 \wr@@x{\string\labeldef{S}\space{#2}\space{#1}}\fi}

\newif\ifd@bn@\d@bn@true
\def\Section{\gad\Sno\ifd@bn@\Fno\z@\Lno\z@\fi\@fn@xt[\s@ct@@\s@ct@}
\def\section{\gad\Sno\ifd@bn@\Fno\z@\Lno\z@\fi\@fn@xt[\s@ct@e\s@cte}
 
\def\subsection{\@fn@xt*\subs@ct@\subs@ct}
\def\subs@ct#1{{\s@bs@ct\empty{#1}}\wr@@c{\string\subcd{#1}{\the\pageno}}}
\def\subs@ct@*#1{\vsk->\nobreak
 {\s@bs@ct\empty{#1}}\wr@@c{\string\subcd{#1}{\the\pageno}}}

\def\l@b@l#1#2{\def\n@@{\csname #2no\endcsname}%
 \if*#1\gad\n@@ \@ft@\xdef\csname @#1@#2@\endcsname{\l@f@rm}\else\def\t@st{#1}%
 \ifx\t@st\empty\gad\n@@ \@ft@\xdef\csname @#1@#2@\endcsname{\l@f@rm}%
 \else\@ft@\ifx\csname @#1@#2@mark\endcsname\relax\gad\n@@
 \@ft@\xdef\csname @#1@#2@\endcsname{\l@f@rm}%
 \@ft@\gdef\csname @#1@#2@mark\endcsname{}%
 \wr@@x{\string\labeldef{#2}\space{\?#1@#2?}\space\ifnum\n@@<10 \space\fi{#1}}%
 \fi\fi\fi}
\def\labeldef#1#2#3{\dff\?#3@#1?{#2}}
\def\Labeldef#1#2#3{\dff\?#3@#1?{#2}\@ft@\gdef\csname @#3@#1@mark\endcsname{}}

\def\l@bel#1#2{\l@b@l{#1}{#2}\?#1@#2?}

\newcount\c@cite
\def\?#1?{\csname @#1@\endcsname}
\def\[{\@fn@xt:\c@t@sect\c@t@}
\def\c@t@#1]{{\c@cite\z@\@fwd@@{\?#1@L?}{\adv\c@cite1}{}%
 \@fwd@@{\?#1@F?}{\adv\c@cite1}{}\@fwd@@{\?#1?}{\adv\c@cite1}{}%
 \relax\ifnum\c@cite=\z@{\bf ???}\wrs@x{No label [#1]}\else
 \ifnum\c@cite=1\let\@@PS\relax\let\@@@\relax\else\let\@@PS\underbar
 \def\@@@{{\rm<}}\fi\@@PS{\?#1?\@@@\?#1@L?\@@@\?#1@F?}\fi}}
\def\(#1){{\rm(\c@t@#1])}}
\def\c@t@s@ct#1{\@fwd@@{\?#1@S?}{\?#1@S?\relax}%
 {{\bf ???}\wrs@x{No section label {#1}}}}
\def\c@t@sect:#1]{\c@t@s@ct{#1}} \let\SNo\c@t@s@ct

\newdimen\l@ftcd \newdimen\r@ghtcd \let\nlc\relax
\newcount\Cdentry \newcount\cdentry \let\prentcd\relax \let\postentcd\relax

\def\d@tt@d{\leaders\hbox to 1em{\kern.1em.\hfil}\hfill}
\def\entcd#1#2#3{\gad\cdentry\prentcd\item{\l@bcdf@nt#1}{\entcdf@nt#2}\alb
 \kern.9em\hbox{}\kern-.9em\d@tt@d\kern-.36em{\p@g@cdf@nt#3}\kern-\r@ghtcd
 \hbox{}\postentcd\par}
\def\Entcd#1#2#3{\gad\Cdentry\global\cdentry\z@
 \def\l@@b{\@fwd@@{#1}{\c@l@b{#1}}{}}\vsk.2>\entcd{\l@@b}{#2}{#3}}
\def\subcd#1#2{{\adv\leftskip.333em\entcd{}{\s@bcdf@nt#1}{#2}}}
\def\Refcd#1#2{\def\t@@st{#1}\ifx\t@@st\empty\ifx\r@fl@ne\empty\relax\else
 \R@fcd{\r@fl@ne}{#2}\fi\else\R@fcd{#1}{#2}\fi}
\def\R@fcd#1#2{\sk@@p{.6\bls}\entcd{}{\hglue-\l@ftcd\R@fcdf@nt #1}{#2}}
\def\Refline{\def\r@fl@ne} \def\Refempty{\let\r@fl@ne\empty}
\def\Appencd{\par\adv\leftskip-\l@ftcd\adv\rightskip-\r@ghtcd\@ppl@ne
 \adv\leftskip\l@ftcd\adv\rightskip\r@ghtcd\let\Entcd\Appcd}
\def\appline{\def\@ppl@ne} \def\Appempty{\let\@ppl@ne\empty}
\def\Appline#1{\def\@ppl@ne{\s@bs@ct{}{#1}}}
\def\Leftcd#1{\adv\leftskip-\l@ftcd\s@twd@\l@ftcd{\c@l@b{#1}\enspace}
 \adv\leftskip\l@ftcd}
\def\Rightcd#1{\adv\rightskip-\r@ghtcd\s@twd@\r@ghtcd{#1\enspace}
 \adv\rightskip\r@ghtcd}
\def\C@nt{Contents} \def\Ap@s{Appendices} \def\R@fcs{References}
\def\contents{\@fn@xt*\cont@@\cont@}
\def\cont@{\@fn@xt[\cnt@{\cnt@[\C@nt]}}
\def\cont@@*{\@fn@xt[\cnt@@{\cnt@@[\C@nt]}}
\def\cnt@[#1]{\c@nt@{M}{#1}{44}{\s@bs@ct{}{\@ppl@f@nt\Ap@s}}}
\def\cnt@@[#1]{\c@nt@{M}{#1}{44}{}}
\def\endco{\par\penalty-500\vsk>\vskip\z@\endgroup}
\def\readcd{\@np@t{\jobname.cd}}
\def\Cde{\@fn@xt*\Cde@@\Cde@}
\def\Cde@{\@fn@xt[\Cd@{\Cd@[\C@nt]}}
\def\Cde@@*{\@fn@xt[\Cd@@{\Cd@@[\C@nt]}}
\def\Cd@[#1]{\cnt@[#1]\readcd\endco}
\def\Cd@@[#1]{\cnt@@[#1]\readcd\endco}
\def\contlabeldef{\def\c@l@b}

\long\def\c@nt@#1#2#3#4{\s@twd@\l@ftcd{\c@l@b{#1}\enspace}
 \s@twd@\r@ghtcd{#3\enspace}\adv\r@ghtcd1.333em
 \def\@ppl@ne{#4}\def\r@fl@ne{\R@fcs}\s@ct{}{#2}\B@gr@\parindent\z@\let\nlc\nl
 \let\nl\relax\parskip.2\bls\adv\leftskip\l@ftcd\adv\rightskip\r@ghtcd}

\def\writecd{\immediate\openout\@@cd\jobname.cd \def\wr@@c{\write\@@cd}
 \def\cl@secd{\immediate\write\@@cd{\string\endinput}\immediate\closeout\@@cd}
 \def\closecd{\cl@secd\global\let\cl@secd\relax}}
\let\cl@secd\relax \def\wr@@c#1{} \let\closecd\relax

\def\dff{\@ft@\d@f} \def\d@f{\@ft@\def}
\def\edff{\@ft@\ed@f} \def\ed@f{\@ft@\edef}
\def\gdff{\@ft@\gd@f} \def\gd@f{\@ft@\gdef}
\def\defi#1#2{\def#1{#2}\wr@@x{\string\def\string#1{#2}}}

\def\qed{\hbox{}\nobreak\hfill\nobreak{\m@th$\,\square$}}
\def\back#1 {\strut\kern-.33em #1\enspace\ignore} 

\def\hcor#1{\advance\hoffset by #1}
\def\vcor#1{\advance\voffset by #1}
\let\bls\baselineskip \let\ignore\ignorespaces
\ifx\ic@\undefined \let\ic@\/\fi
\def\vsk#1>{\vskip#1\bls} \let\adv\advance
\def\vv#1>{\vadjust{\vsk#1>}\ignore}
\def\vvn#1>{\vadjust{\nobreak\vsk#1>\nobreak}\ignore}
\def\vvv#1>{\vskip\z@\vsk#1>\nt\ignore}
\def\vvgood{\vadjust{\penalty-500}}
\def\nngood{\noalign{\penalty-500}}

\def\Goodbreak{\par\penalty-\@m}
\def\wgood#1>{\vv#1>\vvgood\vv-#1>}
\def\wwgood#1:#2>{\vv#1>\vvgood\vv#2>}
\def\mmgood#1:#2>{\cnn#1>\nngood\cnn#2>}
\def\goodsk#1:#2>{\vsk#1>\goodbreak\vsk#2>\vsk0>}
\def\ragood{\vadjust{\vskip\z@ plus 12pt}\vvgood}

 \def\setparindent{\edef\Parindent{\the\parindent}}
\def\Type{\vsk.5>\bgroup\parindent\z@\tt\rightskip\z@ plus1em minus1em%
 \spaceskip.3333em \xspaceskip.5em\relax}
\def\endType{\vsk.5>\egroup\nt} 

 \let\Tilde\widetilde \let\dollar\$ \let\ampersand\&
\let\sss\scriptscriptstyle  
\let\vp\vphantom \let\hp\hphantom \let\nt\noindent
\let\cline\centerline \let\lline\leftline \let\rline\rightline
\def\nn#1>{\noalign{\vskip#1\p@@}} \def\NN#1>{\openup#1\p@@}
\def\cnn#1>{\noalign{\vsk#1>}}
 
\let\Lim\lim \def\lim{\Lim\limits} \let\Sum\sum \def\sum{\Sum\limits}
 
\let\Prod\prod \def\prod{\Prod\limits} \let\Int\int \def\int{\Int\limits}

\def\tprod{\mathop{\tsize\Prod}\limits} \def\&{.\kern.1em}
\def\nl{\leavevmode\hfill\break} \def\~{\leavevmode\@fn@xt~\m@n@s\@md@@sh}
\def\@md@@sh{\@fn@xt-\d@@sh\@md@sh} \def\@md@sh{\raise.13ex\hbox{--}}
\def\m@n@s~{\raise.15ex\mbox{-}} \def\d@@sh-{\raise.15ex\hbox{-}}

\let\procent\% \def\%#1{\ifmmode\mathop{#1}\limits\else\procent#1\fi}
\let\@ml@t\" \def\"#1{\ifmmode ^{(#1)}\else\@ml@t#1\fi}
\let\@c@t@\' \def\'#1{\ifmmode _{(#1)}\else\@c@t@#1\fi}
\let\colon\: \def\:{^{\vp{\topsmash|}}} 

\let\texspace\ \def\ {\ifmmode\alb\fi\texspace} \def\.{\d@t\ignore}

\newif\ifNewskips

\def\Newskips{\global\Newskipstrue
 \gdef\>{\RIfM@\mskip.666667\thinmuskip\relax\else\kern.111111em\fi}
 \gdef\}{\RIfM@\mskip-.666667\thinmuskip\relax\else\kern-.111111em\fi}
 \gdef\){\RIfM@\mskip.333333\thinmuskip\relax\else\kern.0555556em\fi}
 \gdef\]{\RIfM@\mskip-.333333\thinmuskip\relax\else\kern-.0555556em\fi}}
\def\d@t{\ifNewskips.\hskip.3em\else\def\d@t{.\ }\fi} \def\.{\d@t\ignore}
\Newskips

\let\n@wp@ge\newpage \def\newpage{\endgraf\n@wp@ge}
\let\=\m@th \def\mbox#1{\hbox{\m@th$#1$}}
\def\mtext#1{\text{\m@th$#1$}} \def\^#1{\text{\m@th#1}}
\def\Line#1{\kern-.5\hsize\line{\m@th$\dsize#1$}\kern-.5\hsize}
\def\Lline#1{\kern-.5\hsize\lline{\m@th$\dsize#1$}\kern-.5\hsize}
\def\Cline#1{\kern-.5\hsize\cline{\m@th$\dsize#1$}\kern-.5\hsize}
\def\Rline#1{\kern-.5\hsize\rline{\m@th$\dsize#1$}\kern-.5\hsize}

\def\Ll@p#1{\llap{\m@th$#1$}} \def\Rl@p#1{\rlap{\m@th$#1$}}
 \def\Cl@p#1{\llap{\m@th$#1$\hss}}
\def\Llap#1{\mathchoice{\Ll@p{\dsize#1}}{\Ll@p{\tsize#1}}{\Ll@p{\ssize#1}}%
 {\Ll@p{\sss#1}}}
\def\Clap#1{\mathchoice{\Cl@p{\dsize#1}}{\Cl@p{\tsize#1}}{\Cl@p{\ssize#1}}%
 {\Cl@p{\sss#1}}}
\def\Rlap#1{\mathchoice{\Rl@p{\dsize#1}}{\Rl@p{\tsize#1}}{\Rl@p{\ssize#1}}%
 {\Rl@p{\sss#1}}}
 
\def\LRtph#1#2{\setbox\z@\hbox{#1}\dimen\z@\wd\z@\hbox{\hbox to\dimen\z@{#2}}}
\def\LRph#1#2{\LRtph{\m@th$#1$}{\m@th$#2$}}

\def\CCph#1#2{\LRph{#1}{\hss#2\hss}}
 \def\RRph#1#2{\LRph{#1}{#2\hss}}
\def\Cph#1#2{\mathchoice{\CCph{\dsize#1}{\dsize#2}}{\CCph{\tsize#1}{\tsize#2}}
 {\CCph{\ssize#1}{\ssize#2}}{\CCph{\sss#1}{\sss#2}}}

\def\Rph#1#2{\mathchoice{\RRph{\dsize#1}{\dsize#2}}{\RRph{\tsize#1}{\tsize#2}}
 {\RRph{\ssize#1}{\ssize#2}}{\RRph{\sss#1}{\sss#2}}}
\def\Lto#1{\setbox\z@\mbox{\tsize{#1}}%
 \mathrel{\mathop{\hbox to\wd\z@{\rightarrowfill}}\limits#1}}
\def\Lgets#1{\setbox\z@\mbox{\tsize{#1}}%
 \mathrel{\mathop{\hbox to\wd\z@{\leftarrowfill}}\limits#1}}

\let\alb\allowbreak

 \let\x\times \let\ox\otimes 
  \let\tabs\+
\let\le\leqslant \let\ge\geqslant
\let\der\partial \let\8\infty \let\*\star

 \def\vert{\ |\ }

\let\lb\lbrace \let\rb\rbrace

\def\lsym#1{#1\alb\ldots\relax#1\alb}
\def\lc{\lsym,}  \def\lx{\lsym\x} \def\lox{\lsym\ox}

\def\Re{\mathop{\roman{Re}\>}} \def\Im{\mathop{\roman{Im}\>}}

\def\Sym{\mathop{\roman{Sym}\)}\nolimits}

\def\Sing{\mathop{\roman{Sing}\>}}

 \def\1{^{-1}} \let\underscore\_ \def\_#1{_{\Rlap{#1}}}
\def\vst#1{{\lower1.9\p@@\mbox{\bigr|_{\raise.5\p@@\mbox{\ssize#1}}}}}
\def\vrp#1:#2>{{\vrule height#1 depth#2 width\z@}}
\def\vru#1>{\vrp#1:\z@>} \def\vrd#1>{\vrp\z@:#1>}
\def\qqq{\qquad\quad} 
\def\sscr#1{\raise.3ex\mbox{\sss#1}} \def\@@PS{\bold{OOPS!!!}}
\def\ono{\bigl(1+o(1)\bigr)} 

\def\intcl{\mathop
 {\Rlap{\raise.3ex\mbox{\kern.12em\curvearrowleft}}\int}\limits}
\def\intcr{\mathop
 {\Rlap{\raise.3ex\mbox{\kern.24em\curvearrowright}}\int}\limits}

\def\pms{\raise.25ex\mbox{\ssize\pm}\>}
\def\mps{\raise.25ex\mbox{\ssize\mp}\>}

\let\al\alpha
\let\bt\beta
\let\gm\gamma \let\Gm\Gamma 
 \let\Dl\Delta 
 \let\eps\varepsilon \let\epsilon\eps

\let\la\lambda 

\let\si\sigma 
 
 \let\phi\varphi

\let\om\omega \let\Om\Omega 

\def\C{\Bbb C}
\def\R{\Bbb R}
\def\Z{\Bbb Z}

\def\difl/{differential} \def\dif/{difference}
\def\cf.{cf.\ \ignore} \def\Cf.{Cf.\ \ignore}
\def\egv/{eigenvector} \def\eva/{eigenvalue} \def\eq/{equation}
\def\lhs/{the left hand side} \def\rhs/{the right hand side}
\def\Lhs/{The left hand side} \def\Rhs/{The right hand side}
\def\gby/{generated by} \def\wrt/{with respect to} \def\st/{such that}
\def\resp/{respectively} \def\off/{offdiagonal} \def\wt/{weight}
\def\pol/{polynomial} \def\rat/{rational} \def\tri/{trigonometric}
\def\fn/{function} \def\var/{variable} \def\raf/{\rat/ \fn/}
\def\inv/{invariant} \def\hol/{holomorphic} \def\hof/{\hol/ \fn/}
\def\mer/{meromorphic} \def\mef/{\mer/ \fn/} \def\mult/{multiplicity}
\def\sym/{symmetric} \def\perm/{permutation}
\def\rep/{representation} \def\irr/{irreducible} \def\irrep/{\irr/ \rep/}
\def\hom/{homomorphism} \def\aut/{automorphism} \def\iso/{isomorphism}
\def\lex/{lexicographical} \def\as/{asymptotic} \def\asex/{\as/ expansion}
\def\ndeg/{nondegenerate} \def\neib/{neighbourhood} \def\deq/{\dif/ \eq/}
\def\hw/{highest \wt/} \def\gv/{generating vector} \def\eqv/{equivalent}
\def\msd/{method of steepest descend} \def\pd/{pairwise distinct}
\def\wlg/{without loss of generality} \def\Wlg/{Without loss of generality}
\def\onedim/{one-dim\-en\-sion\-al} \def\fd/{fin\-ite-dim\-en\-sion\-al}
\def\qcl/{quasiclassical} \def\hwv/{\hw/ vector}
\def\hgeom/{hyper\-geo\-met\-ric} \def\hint/{\hgeom/ integral}
\def\hwm/{\hw/ module} \def\emod/{evaluation module} \def\Vmod/{Verma module}
\def\symg/{\sym/ group} \def\sol/{solution} \def\eval/{evaluation}
\def\anf/{analytic \fn/} \def\anco/{analytic continuation}
\def\qg/{quantum group} \def\qaff/{quantum affine algebra}

\hyphenation{ortho-gon-al}

\def\Rm/{\^{$R$-}matrix} \def\Rms/{\^{$R$-}matrices} \def\YB/{Yang-Baxter \eq/}
\def\Ba/{Bethe ansatz} \def\Bv/{Bethe vector} \def\Bae/{\Ba/ \eq/}
\def\KZv/{Knizh\-nik-Zamo\-lod\-chi\-kov} \def\KZvB/{\KZv/-Bernard}
\def\KZ/{{\sl KZ\/}} \def\qKZ/{{\sl qKZ\/}}
\def\KZB/{{\sl KZB\/}} \def\qKZB/{{\sl qKZB\/}}
\def\qKZo/{\qKZ/ operator} \def\qKZc/{\qKZ/ connection}
\def\KZe/{\KZ/ \eq/} \def\qKZe/{\qKZ/ \eq/} \def\qKZBe/{\qKZB/ \eq/}
\def\XXX/{{\sl XXX\/}} \def\XXZ/{{\sl XXZ\/}} \def\XYZ/{{\sl XYZ\/}}

\def\h@ph{\discretionary{}{}{-}} \def\$#1$-{\,\^{$#1$}\h@ph}

\def\TFT/{Research Insitute for Theoretical Physics}
\def\HY/{University of Helsinki} \def\AoF/{the Academy of Finland}
\def\CNRS/{Supported in part by MAE\~MICECO\~CNRS Fellowship}
\def\LPT/{Laboratoire de Physique Th\'eorique ENSLAPP}
\def\ENSLyon/{\'Ecole Normale Sup\'erieure de Lyon}
\def\LPTaddr/{46, All\'ee d'Italie, 69364 Lyon Cedex 07, France}
\def\enslapp/{URA 14\~36 du CNRS, associ\'ee \`a l'E.N.S.\ de Lyon,
au LAPP d'Annecy et \`a l'Universit\`e de Savoie}
\def\ensemail/{vtarasov\@ enslapp.ens-lyon.fr}
\def\DMS/{Department of Mathematics, Faculty of Science}
\def\DMO/{\DMS/, Osaka University}
\def\DMOaddr/{Toyonaka, Osaka 560, Japan}
\def\dmoemail/{vt\@ math.sci.osaka-u.ac.jp}
\def\MPI/{Max\)-Planck\)-Institut} \def\MPIM/{\MPI/ f\"ur Mathematik}
\def\MPIMaddr/{P\]\&O.\ Box 7280, D\~-\]53072 \,Bonn, Germany}
\def\mpimemail/{tarasov\@ mpim-bonn.mpg.de}
\def\SPb/{St\&Peters\-burg}
\def\home/{\SPb/ Branch of Steklov Mathematical Institute}
\def\homeaddr/{Fontanka 27, \SPb/ \,191011, Russia}
\def\homemail/{vt\@ pdmi.ras.ru}
\def\absence/{On leave of absence from \home/}
\def\support/{Supported in part by}
\def\UNC/{Department of Mathematics, University of North Carolina}
\def\ChH/{Chapel Hill} \def\UNCaddr/{\ChH/, NC 27599, USA}
\def\avemail/{anv\@ email.unc.edu}	
\def\grant/{NSF grant DMS\~9501290}	
\def\Grant/{\support/ \grant/}

\def\Aomoto/{K\&Aomoto}
\def\Cher/{I\&Che\-red\-nik}
\def\Dri/{V\]\&G\&Drin\-feld}
\def\Fadd/{L\&D\&Fad\-deev}
\def\Feld/{G\&Felder}
\def\Fre/{I\&B\&Fren\-kel}
\def\Etingof/{P\]\&Etingof}
\def\Gustaf/{R\&A\&Gustafson}
\def\Izergin/{A\&G\&Izergin}
\def\Jimbo/{M\&Jimbo}
\def\Kazh/{D\&Kazhdan}
\def\Kor/{V\]\&E\&Kore\-pin}
\def\Kulish/{P\]\&P\]\&Ku\-lish}
\def\Lusz/{G\&Lusztig}
\def\Miwa/{T\]\&Miwa}
\def\MN/{M\&Naza\-rov}
\def\Reshet/{N\&Reshe\-ti\-khin} \def\Reshy/{N\&\]Yu\&Reshe\-ti\-khin}
\def\SchV/{V\]\&\]V\]\&Schecht\-man} \def\Sch/{V\]\&Schecht\-man}
\def\Skl/{E\&K\&Sklya\-nin}
\def\Smirnov/{F\]\&A\&Smir\-nov}
\def\Takh/{L\&A\&Takh\-tajan}
\def\VT/{V\]\&Ta\-ra\-sov} \def\VoT/{V\]\&O\&Ta\-ra\-sov}
\def\Varch/{A\&\]Var\-chenko} \def\Varn/{A\&N\&\]Var\-chenko}
\def\Zhel/{D\&P\]\&Zhe\-lo\-ben\-ko}

\def\AiA/{Al\-geb\-ra i Ana\-liz}
\def\DAN/{Do\-kla\-dy AN SSSR}
\def\FAA/{Funk\.Ana\-liz i ego pril.}
\def\Izv/{Iz\-ves\-tiya AN SSSR, ser\.Ma\-tem.}
\def\TMF/{Teo\-ret\.Ma\-tem\.Fi\-zi\-ka}
\def\UMN/{Uspehi Matem.\ Nauk}

\def\AMS/{Amer\.Math\.Society}
\def\AMSa/{AMS \publaddr Providence RI}
\def\AMST/{\AMS/ Transl.,\ Ser\&\)2}
\def\AMSTr/{\AMS/ Transl.,} \def\Ser2{Ser\&\)2}
\def\Astq/{Ast\'erisque}
\def\ContM/{Contemp\.Math.}
\def\CMP/{Comm\.Math\.Phys.}
\def\DMJ/{Duke\.Math\.J.}
\def\Inv/{Invent\.Math.} 
\def\IMRN/{Int\.Math\.Res.\ Notices}
\def\JMP/{J\.Math\.Phys.}
\def\JPA/{J\.Phys.\ A}
\def\JSM/{J\.Soviet Math.}
\def\LMJ/{Leningrad Math.\ J.}
\def\LpMJ/{\SPb/ Math.\ J.}
\def\LMP/{Lett\.Math\.Phys.}
\def\NMJ/{Nagoya Math\.J.}
\def\Nucl/{Nucl\.Phys.\ B}
\def\OJM/{Osaka J\.Math.}
\def\RIMS/{Publ\.RIMS, Kyoto Univ.}
\def\SIAM/{SIAM J\.Math\.Anal.}
\def\SMNS/{Selecta Math., New Series}
\def\TMP/{Theor\.Math\.Phys.}
\def\ZNS/{Zap\. nauch\. semin. LOMI}

\def\ASMP/{Advanced Series in Math.\ Phys.{}}

\def\Birk/{Birkh\"auser}
\def\CUP/{Cambridge University Press} \def\CUPa/{\CUP/ \publaddr Cambridge}
\def\Spri/{Springer\)-Verlag} \def\Spria/{\Spri/ \publaddr Berlin}
\def\WS/{World Scientific} \def\WSa/{\WS/ \publaddr Singapore}

\newbox\lefthbox \newbox\righthbox

\let\sectsep. \let\labelsep. \let\contsep. \let\labelspace\relax
\let\sectpre\relax \let\contpre\relax
\def\sf@rm{\the\Sno} \def\sf@rm@{\sectpre\sf@rm\sectsep}
\def\c@l@b#1{\contpre#1\contsep}
\def\l@f@rm{\ifd@bn@\sf@rm\labelsep\fi\labelspace\the\n@@}

\def\sectformdef{\def\sf@rm}

\let\DoubleNum\d@bn@true \let\SingleNum\d@bn@false

\def\NoNewNum{\let\writeldf\relax\def\l@b@l##1##2{\if*##1%
 \@ft@\xdef\csname @##1@##2@\endcsname{\mbox{*{*}*}}\fi}}
\def\NoNewTime{\def\todaydef##1{\def\today{##1}}
 \def\nowtimedef##1{\def\nowtime{##1}}}
\def\NoInput{\let\Input\input\let\writeldf\relax}
\def\Fixed{\NoNewTime\NoInput}

\newbox\dtlb@x
\def\DateTimeLabel{\global\setbox\dtlb@x\vbox to\z@{\ifMag\eightpoint\else
 \ninepoint\fi\sl\vss\rline\today\rline\nowtime}
 \global\headline{\hfil\box\dtlb@x}}

\def\sectfont#1{\def\s@cf@nt{#1}} \sectfont\bf
\def\subsectfont#1{\def\s@bf@nt{#1}} \subsectfont\it
\def\Entcdfont#1{\def\entcdf@nt{#1}} \Entcdfont\relax
\def\labelcdfont#1{\def\l@bcdf@nt{#1}} \labelcdfont\relax
\def\pagecdfont#1{\def\p@g@cdf@nt{#1}} \pagecdfont\relax
\def\subcdfont#1{\def\s@bcdf@nt{#1}} \subcdfont\it
\def\applefont#1{\def\@ppl@f@nt{#1}} \applefont\bf
\def\Refcdfont#1{\def\R@fcdf@nt{#1}} \Refcdfont\bf

\def\reffont#1{\def\r@ff@nt{#1}} \reffont\rm
\def\keyfont#1{\def\k@yf@nt{#1}} \keyfont\rm
\def\paperfont#1{\def\p@p@rf@nt{#1}} \paperfont\it
\def\bookfont#1{\def\b@@kf@nt{#1}} \bookfont\it
\def\volfont#1{\def\v@lf@nt{#1}} \volfont\bf
\def\issuefont#1{\def\iss@f@nt{#1}} \issuefont{no\p@@nt}

\def\adjustmid#1{\kern-#1\p@\alb\hskip#1\p@\relax}
\def\adjustend#1{\adjustnext{\kern-#1\p@\alb\hskip#1\p@}}

\newif\ifcd 

\tenpoint

\Fixed

\Magset
\ifx\USversion\relax\UStrue\fi

\PaperA4
\ifMag
\ifUS
\PaperUS
\hoffset\z@
\voffset\z@
\fi
\else
\USfalse
\fi

\font\symf=cmsy10 scaled\magstep5

\def\ub{\bar u}
\def\vb{\bar v}

\def\Sb{\bold S}

\def\Cc{\Cal C}

\def\gsl{\frak{sl}}

\def\Gt{{\>\Tilde{\}G\>}\}}}
\def\hti{\Rlap{\>\Tilde{\}\phantom h\)}\]}h}
\def\Jt{\Rlap{\,\Tilde{\!\phantom J\,\)}\]\!}J}
\def\Rp{\R_{\ge 0}}

\def\Syml#1{\Rlap{\Sym\limits_{\Rph{\tsize\Sym}{\)#1}}}\hp{\ssize\)#1}}

\def\veeq{\raise.36ex\hbox{\ninepoint\=$\ssize{\backslash}\!/\!\]/$}}

\def\wedq{\raise.4ex\hbox{\ninepoint\=$\ssize/\!{\backslash}\!\]{\backslash}$}}

\def\semibold#1{\setbox\p@b@ld\hbox{#1}\kern-.006em\copy\p@b@ld\kern-\wd\p@b@ld
 \kern.013em\copy\p@b@ld\kern-\wd\p@b@ld\kern-.007em\raise.013em\box\p@b@ld}

\def\halfbold#1{\setbox\p@b@ld\hbox{#1}\kern-.005em\copy\p@b@ld\kern-\wd\p@b@ld
 \kern.0106em\copy\p@b@ld\kern-\wd\p@b@ld\kern-.005em\raise.011em\box\p@b@ld}

\def\vveq{\raise.36ex\mbox{\semibold{\ninepoint\=$\ssize{\backslash}$}
 \sssize\!\halfbold{\ninepoint\=$\ssize/\!\]/$}}}

\def\wwdq{\raise.4ex\mbox{\halfbold{\ninepoint\=$\ssize/$}\sssize\!
 \semibold{\ninepoint\=$\ssize{\backslash}\!\]{\backslash}$}}}

\def\roop{\vcenter{\hbox{\symf\char26\kern-10.5\p@
 \vbox{\parindent\z@ \hrule height 1\p@ depth \z@ width .15\hsize
 \hrule height 6.02\p@ depth \z@ width \z@
 \hrule width .15\hsize
 \hrule height 5.98\p@ depth .04\p@ width \z@
 \hrule height -.04\p@ depth 1.04\p@ width .15\hsize}\kern-.4em
 \smash{\lower3.5\p@\vbox{\m@th\twelvepoint\mbox{\hp{\to}\!\!\}{\leftarrow}}
 \vskip-.9\p@\mbox{{\to}\!\!\}\hp{\leftarrow}}}}\kern-.4em
  \vbox{\parindent\z@ \hrule height 1\p@ depth \z@ width .15\hsize
 \hrule height 6.02\p@ depth \z@ width \z@
 \hrule width .15\hsize
 \hrule height 5.98\p@ depth .04\p@ width \z@
 \hrule height -.04\p@ depth 1.04\p@ width .15\hsize}
 \kern-.1em\smash{\raise3.7\p@\vbox{\m@th\mbox{\to}}}}}}

\def\dkt{d^{\)k}\]t}
\def\dkts{d^{\)k_1}\]t\;d^{\)k_2}\}s}

\def\msm{{\sssize M}}

\def\Dlk{\Dl^{\]k}_{\vp1}}
\def\Dlkot{\Dl^{\]k_1\],\)k_2}_{\vp1}}
\def\DlkZ{\mathchoice{\Dl^{\]k}_{\)\Z,\gm}}{\Dl^{\]k}_{\)\Z,\gm}}
 {\Dl^{\]k}_{\)\Z\),\]\gm}}{\@@PS}}
\def\DlkotZ{\mathchoice{\Dl^{\]k_1\],\)k_2}_{\)\Z\),\gm}}
 {\Dl^{\]k_1\],\)k_2}_{\)\Z,\gm}}{\Dl^{\]k_1\],k_2}_{\)\Z,\]\gm}}{\@@PS}}

\def\Ckot{C^{\)k_1\],\)k_2}_\gm}
\def\Dkot{D^{\)k_1\],\)k_2}_{\]\msm}}
\def\Rkot{\R^{\)k_1+\)k_2}_{\vp1}}
\def\Xkot{X^{k_1\],\)k_2}_{\msm\],\gm}}

\def\zn{z_1\lc z_n}
\def\Vox{V_1\lox V_n}
\def\tk{t_1\lc t_k}
\def\skt{s_1\lc s_{k_2}}
\def\tko{t_1\lc t_{k_1}}
\def\uko{u_1\lc u_{k_1}}
\def\vkt{v_1\lc v_{k_2}}

\def\prak{\prod_{a=1}^k}
\def\prabk{\prod_{1\le a<b\le k}}
\def\prako{\prod_{a=1}^{k_1}}
\def\prbkt{\prod_{b=1}^{k_2}}
\def\prabko{\prod_{1\le a<b\le k_1}}
\def\prabkt{\prod_{1\le a<b\le k_2}}
\def\prbakt{\prod_{1\le b<a\le k_2}}

\let\goodbm\relax  \let\mmgood\relax 
  \let\uugood\relax 
 \def\vvm#1>{\ignore} \def\vvnm#1>{\ignore} \def\cnnm#1>{}
\def\cnnu#1>{} \def\vvu#1>{\ignore} \def\vvnu#1>{\ignore} 

\ifMag \ifUS   \let\uugood\vvgood
  \let\vvu\vv \let\vvnu\vvn \let\cnnu\cnn  \else
 \let\goodbm\goodbreak  \let\mmgood\vvgood \let\cnnm\cnn
  \let\vvm\vv \let\vvnm\vvn  \fi
 \let\goodbreak\relax  \let\vvgood\relax
  \def\nngood{}  \fi

\def\wwmgood#1:#2>{\ifMag\vv#1>\mmgood\vv#2>\vv0>\fi}
\def\vskmgood#1:#2>{\ifMag\vsk#1>\goodbm\vsk#2>\vsk0>\fi}
\def\vskm#1:#2>{\ifMag\vsk#1>\else\vsk#2>\fi}
\def\vvmm#1:#2>{\ifMag\vv#1>\else\vv#2>\fi}
\def\vvnn#1:#2>{\ifMag\vvn#1>\else\vvn#2>\fi}
\def\nnm#1:#2>{\ifMag\nn#1>\else\nn#2>\fi}
\def\kerm#1:#2>{\ifMag\kern#1em\else\kern#2em\fi}

\def\volume{vol\point}
\sectfont{\elevenpoint\bf}

\csname selberg.def\endcsname

\labeldef{F} {1\labelsep \labelspace 1}  {selb}
\labeldef{F} {1\labelsep \labelspace 2}  {exp}
\labeldef{F} {1\labelsep \labelspace 3}  {dexp}
\labeldef{F} {1\labelsep \labelspace 4}  {zkg}

\labeldef{F} {2\labelsep \labelspace 1}  {phi}
\labeldef{F} {2\labelsep \labelspace 2}  {wtf}
\labeldef{F} {2\labelsep \labelspace 3}  {Fun}
\labeldef{L} {2\labelsep \labelspace 1}  {Flim}
\labeldef{F} {2\labelsep \labelspace 4}  {dlk12}
\labeldef{L} {2\labelsep \labelspace 2}  {Fval}
\labeldef{L} {2\labelsep \labelspace 3}  {sl3dexp}
\labeldef{F} {2\labelsep \labelspace 5}  {dexp3}

\labeldef{F} {3\labelsep \labelspace 1}  {wtg}
\labeldef{L} {3\labelsep \labelspace 1}  {sl3exp}
\labeldef{F} {3\labelsep \labelspace 2}  {exp3}
\labeldef{L} {3\labelsep \labelspace 2}  {sl3selb}
\labeldef{F} {3\labelsep \labelspace 3}  {selb3}
\labeldef{L} {3\labelsep \labelspace 3}  {sl3selb0}
\labeldef{F} {3\labelsep \labelspace 4}  {selb30}
\labeldef{L} {3\labelsep \labelspace 4}  {chains}

\labeldef{L} {5\labelsep \labelspace 1}  {more}
\labeldef{F} {5\labelsep \labelspace 1}  {JJl}
\labeldef{F} {5\labelsep \labelspace 2}  {JJ0}
\labeldef{L} {5\labelsep \labelspace 2}  {JJJ}
\labeldef{F} {5\labelsep \labelspace 3}  {J1}
\labeldef{F} {5\labelsep \labelspace 4}  {Jt1}
\labeldef{F} {5\labelsep \labelspace 5}  {J0k}

\document

\center
\hrule height 0pt
\vsk-.2>

{\twelvepoint\bf \bls1.2\bls
Selberg Type Integrals Associated with $\gsl_3$
\par}
\vsk1.5>
\=
\VT/$^{\,\star}$ \ and \ \Varch/$^{\,*}$
\vsk1.5>
{\it $^\star$\home/\\ \homeaddr/
\vsk.3>
$^*$\UNC/\\ \UNCaddr/
\vsk1.5>
\sl February \,2003}
\endcenter

\ftext{\=\bls11pt
$\]^\star\)$\support/ RFFI grant 02\)\~\)01\~\)00085a
\>and \,CRDF grant RM1\~\)2334\)\~MO\)\~\)02\vv.06>\nl
\hp{$^*$}{\tenpoint\sl E-mail\/{\rm:} \homemail/}\vv.1>\nl
${\]^*\)}$\support/ NSF grant DMS\)\~\)9801582\vv.06>\nl
\vv-1.2>
\hp{$^*$}{\tenpoint\sl E-mail\/{\rm:} \avemail/}}

\vskm1.4:1.6>

{\ifMag\ninepoint\fi
\Abstract
We present several formulae for Selberg type integrals
associated with the Lie algebra $\gsl_3$.
\endAbs}

\vskm1.2:1.4>
\vsk0>

\Sno 0

\Section{Introduction}
The Selberg integral is the integral
$$
\int_{\Dlk[\)0,1]}
\prak\,t_a^{\>\al-1}\)(1-t_a)^{\)\bt-1}\!\!\prabk\!\}(t_a\]-t_b)^{\)2\gm}\>\dkt
\vv-.2>
\Tag{selb}
$$
where $\Dlk[x\),y]\)=\)\lb\)t\in\R^k\vert x\le t_k\lsym\le t_1\le y\)\rb$.
Integral \(selb) is a generalization of the Euler beta integral.
In 1944 Selberg showed \cite{S} that the integral equals
\vvn.2>
$$
\prod_{j=0}^{k-1}\,{\Gm(\al+j\)\gm)\,\Gm(\)\bt+j\)\gm)\,\Gm(\gm+j\gm)
\over\Gm\bigl(\al+\bt+(2\)k-2-j\))\gm\bigr)\,\Gm(\gm)}\;.
\vv.2>
$$
The Selberg integral is one of the most remarkable \hgeom/ \fn/s with many
applications, see for instance \cite{A1}\), \cite{A2}\), \cite{As}\),
\cite{D}\), \cite{DF1}\), \cite{DF2}\), \cite{FSV}\), \cite{M}\).
Taking a suitable limit of the Selberg integral one gets
the exponential Selberg integral:
$$
\int_{\Dlk[\)0,\)+\8]}
\prak\,e^{-\)t_a}\>t_a^{\>\al-1}\!\!\prabk\!\}(t_a\]-t_b)^{\)2\gm}\>\dkt
\;=\,\prod_{j=0}^{k-1}\,{\Gm(\al+j\gm)\,\Gm(\gm+j\)\gm)\over\Gm(\gm)}\;.
\Tag{exp}
$$
\par
There is also a discrete version of the exponential Selberg integral:
$$
\align
\sum_{u\)\in\)\DlkZ\!}\;
\prak {}& \,z^{u_a}\>{\Gm(u_a\]+\al)\over\Gm(u_a\]+1)}\>
\prabk\!\!{(u_a\]-u_b)\>\Gm(u_a\]-u_b\]+\gm)\over\Gm(u_a\]-u_b\]-\gm+1)}\;={}
\kerm-1:0>
\Tagg{dexp}
\\
{}={}\,\) & \] z^{\)k(k-1)\)\gm/2}\>(1-z)^{-k\al\)-k(k\)-1)\)\gm}\,
\prod_{j=0}^{k-1}\,{\Gm(\al+j\gm)\,\Gm(\gm+j\)\gm)\over\Gm(\gm)}
\kerm-1:0>
\\
\cnn-.2>
\endalign
$$
where ${\DlkZ\]=\)(k-1\),k-2\)\lc 0)\>\gm\>+
\bigl(\)\Dlk[\)0\),+\)\infty]\>\cap\>\Z^k\)\bigr)}$.
Formula \(dexp) can be obtained from the formula for the \$q\)$-Selberg
integral, see \cite{AK}\), by the so-called \rat/ degeneration.
\vsk.1>
The series \(dexp) converges for $|\)z\)|<1$ and any $\al\),\gm$ \st/ all
gamma-\fn/s are well defined. As ${z\to 1}$, formula \(dexp) reproduces
formula \(exp)\).
\vsk.1>
For generic $\al\),\gm$ one can replace the sum over the lattice cone $\DlkZ$
by the sum over the total lattice
$$
\Z^k_\gm\,=\,(k-1\),k-2\)\lc0)\>\gm\>+\>\Z^k,
\Tag{zkg}
$$
since all additional terms appearing in the sum over $\Z^k_\gm$ vanish.
\vsk.2>
All three types of Selberg integrals are related to \rep/ theory of
the Lie algebra $\gsl_2$. Namely, given a tensor product ${V=\)\Vox}$ of
\$\gsl_2\)$-modules one has various linear systems of \difl/ and \deq/s
for a \$V\}$-valued \fn/ $u(\zn)$. The \eq/s are called the \KZv/ (\KZ/\))
and dynamical \eq/s, see~\cite{EFK}\), \cite{EV}\), \cite{FMTV}\), \cite{KZ}\),
\cite{TV2}\), \cite{TV3}\), \cite{V}\). Explicit formulae can be found in
\cite{TV3}\). All three types of Selberg integrals appear as coordinate \fn/s
of \hgeom/ \sol/s of those \eq/s, see~\cite{FMTV}\), \cite{MaV}\),
\cite{TV1}\), \cite{SV}\), \cite{V}\).
\vsk.1>
In particular, integral \(selb) appears as a coordinate \fn/ of
the \hgeom/ \sol/ of the \KZ/ \difl/ \eq/s with values in the space
$$
\Sing(L_{\ell_1}\!\ox L_{\ell_2})\)[\)\ell_1\]+\ell_2-2\)k\)]\,=\,
\bigl\lb\)v\in L_{\ell_1}\!\ox L_{\ell_2}\vert
h\)v=(\ell_1\]+\ell_2-2\)k)\>v\,,\ \,e\)v=0\)\bigr\rb\,,
$$
where $L_\ell$ is the \irr/ \$\gsl_2\)$-module with \hw/ $\ell\in\C$, and
$e\),\)h$ are the standard generators of the Borel subalgebra of $\gsl_2$.
\vsk.1>
Integral \(exp) appears as the coordinate \fn/ of the \hgeom/ \sol/ of
the dynamical \difl/ \eq/ from \cite{FMTV} with values in the \wt/ subspace
$$
L_\ell\)[\)\ell-2\)k\)]\,=\,
\bigl\lb\)v\in L_\ell\vert h\)v=(\ell-2\)k)\>v\)\bigr\rb\,.
$$
Integral \(dexp) appears as the coordinate \fn/ of the \hgeom/ \sol/ of the
dynamical \difl/ \eq/ from \cite{TV3} with values in $L_\ell\)[\)\ell-2\)k\)]$.
\vsk.1>
In all three cases the unknown \fn/ $u(\zn)$ takes values in a \onedim/ space.
It was conjectured in \cite{MV} that if the space of values of a \KZ/ or
a dynamical \eq/ is \onedim/, then coordinates of the \hgeom/ \sol/
can be expressed in terms of elementary \fn/s and gamma-\fn/s.
\vsk.1>
Motivated by this ideology, we give in this paper four versions of new Selberg
type integrals analogous to integrals \(selb)\>--\>\(dexp)\), see formulae
\(dexp3)\), \(exp3)\>--\>\(selb30)\). The new integrals are coordinate \fn/s
of \hgeom/ \sol/s of \KZ/ and dynamical \eq/s associated with the Lie algebra
$\gsl_3$, but we do not discuss this connection in the present paper.
The consideration extends in a rather straightforward way to the $\gsl_n$ case
for $n>3$ giving the corresponding Selberg type integrals. This will be done in
a separate paper.
\vsk.1>
The generalization of the Selberg integral proposed in this paper is
complementary to the definition of Selberg type integrals associated with root
systems \cite{Ma}\). In that definition the root system enters the integrand
via the product of linear \fn/s of integration \var/s, the linear \fn/s
defining the reflection hyperplanes of the root system. From that point of view
the Selberg integral \(selb) is associated with the root system $A_{k-1}$. In
our generalization of the Selberg integral the linear factors of the integrand
do not depend on the root system, but the root system governs the exponents of
linear factors of the integrand, \cf. \(selb3), \(selb30). From our point of
view the Selberg integral, defined for any $k$ by formula \(selb), corresponds
to the root system $A_1$.
\vsk.1>
One can expect that there are Selberg type integrals depending on two root
systems, the first one determining the linear factors of the integrand, like
in \cite{Ma}, and the second one governing the exponents, like in this paper.

\Section{Discrete exponential Selberg integrals associated with $\gsl_3$}

For any \fn/ $f(\tk)$ set
\vvn-.1>
$$
\Syml{\tk}f(\tk)\,=\;
{1\over k!}\,\sum_{\si\)\in\)\Sb_k\!\}}\)f(t_{\si_1}\lc t_{\si_k})\,.
$$
\par
Fix nonnegative integers $k_1\),k_2$ \st/ $k_1\ge k_2$.
Introduce the master \fn/
$$
\align
\Phi({}\!\] &\,\)\uko;\)\vkt;\)z_1,z_2)\,={}
\Tagg{phi}
\\
\nn8>
& {}=\,\prako\,z_1^{u_a}\>{\Gm(u_a\]+\al)\over\Gm(u_a\]+1)}\ \prbkt\,z_2^{v_b}
\,\prako\,\prbkt\,{\Gm(v_b\]-u_a\]-\gm+1)\over\Gm(v_b\]-u_a\]+1)}\;\x{}
\\
\nn7>
&{}\)\x\prabko{(u_a\]-u_b)\>\Gm(u_a\]-u_b\]+\gm)\over\Gm(u_a\]-u_b\]-\gm+1)}
\prabkt{(v_a\]-v_b)\>\Gm(v_a\]-v_b\]+\gm)\over\Gm(v_a\]-v_b\]-\gm+1)}\;,
\kerm-1.6:0>
\\
\cnn.3>
\endalign
$$
and the \wt/ \fn/
\vvn-.3>
$$
\alignat2
w({}\!\] & \,\)\uko;\)\vkt)\,=\,
\Tag{wtf}
\\
\nn6>
& {}=\,\Syml{\uko}\,\Syml{\vkt}
\Bigl(\;\prbkt\,{1\over v_b\]-u_{b\)+k_1\]-k_2}\}-\gm}\,
\prbakt\,{v_b\]-u_{a\)+k_1\]-k_2}\over v_b\]-u_{a\)+k_1\]-k_2}\}-\gm}\;\x{}
\kerm-1.8:0>
\\
\nn6>
&& \Llap{\x\prabko\!\!{u_a\]-u_b\]-\gm\over u_a\]-u_b}
\prabkt\!\!{v_a\]-v_b\]-\gm\over v_a\]-v_b}\,\Bigr)\Rlap{\,.}\kerm-1.8:0>} &
\endalignat
$$
The \wt/ \fn/ is \sym/ in $\uko$ and in $\vkt$ separately, and has
\uugood
at most simple poles located at the hyperplanes $v_b\]-u_a=\gm$. Its numerator
$$
w(\uko;\)\vkt)\,\prako\,\prbkt\,(v_b\]-u_a\]-\gm)
$$
vanishes at every triple intersection of hyperplanes of the form
$\,u_a=\)u_b\]+\gm=\)v_c\>$ or of the form $\,u_a=\)v_b=\)v_c\]-\gm$. Set
\ifMag
$$
\align
& F(\uko;\)\vkt;\)z_1,z_2)\,={}
\Tag{Fun}
\\
\nn4>
&\!\]{}=\,\Phi(\uko;\)\vkt;\)z_1,z_2)\,w(\uko;\)\vkt)\,.
\\
\cnn.2>
\endalign
$$
\else
$$
F(\uko;\)\vkt;\)z_1,z_2)\,=\,\Phi(\uko;\)\vkt;\)z_1,z_2)\,w(\uko;\)\vkt)\,.
\Tag{Fun}
$$
\fi
\Lm{Flim}
Let $(\ub\),\vb)\in\Z^{k_1}_\gm{\x\,}\Z^{k_2}_\gm$, \cf. \(zkg)\).
Let $\al\),\gm$ be generic. Then the \fn/ $F(u\);v\);z_1,z_2)$ has a limit as
$(u\>,v)\to(\ub\>,\vb)$ along any straight line which does not belong to the
singularity hyperplanes of the \fn/ $F$\}, and the limit does not depend on the
direction. By abuse of notation we denote this limit $F(\ub\);\vb\);z_1,z_2)$.
\endpro
Actually, in a neighbourhood of any point
$(\ub\),\vb)\in\Z^{k_1}_\gm{\x\,}\Z^{k_2}_\gm$ the \fn/ $F$ can be written in
the form $F=\)G+P/Q$, where $G\),\)P,\>Q$ are \fn/s, regular at $(\ub\),\vb)$,
and the leading part of the Taylor series expansion of the \fn/ $P$ at
$(\ub\),\vb)$ has greater homogeneous degree than that of the \fn/ $Q$.
Hence, the limit $F(\ub\);\vb\);z_1,z_2)$ equals $G(\ub\);\vb)$.
\vsk.4>
Let $\Dlkot[x,y]$ be a domain in $\Rkot\!$ with coordinates $\tko,\)\skt$
defined by the following inequalities:
\vvn-.1>
$$
\alignat2
x\le{} & t_{k_1}\le t_{k_1-1}\le\;\ldots\;\le t_{k_1-k_2+1}
\le t_{k_1-k_2}\le\;\ldots\;\le t_1 \le y\kerm-2:-1>
\Tag{dlk12}
\\
\nn2>
& \Rph{t_{k_1}}{\wedq}\hp{{}\le{}}\Cph{t_{k_1-1}}{\wedq\)}
\Cph{{}\le\;\ldots\;\le{}}{\cdots}\Rph{s_1}{\wedq} && \kerm2:1>\; .
\\
\nn-3>
& \Cph{t_{k_1}}{s_{k_2}}\le\Rph{t_{k_1-1}}{s_{k_2-1}}\le\;\ldots\;\le s_1\le y
\endalignat
$$
Let \;$\DlkotZ\]=\)(k_1\]-1\)\lc 0\),\)k_2\]-1\),\)\lc 0)\>\gm\>+
\bigl(\)\Dlkot[\)0\),+\)\infty]\>\cap\>\Z^{k_1+\)k_2}_{\vp1}\)\bigr)$.
\Lm{Fval}
Let $(u\),v)\in\Z^{k_1}_\gm{\x\,}\Z^{k_2}_\gm$. Let $\al\),\gm$ be generic.
Then $F(u\);v\);z_1,z_2)\)=\)0$ unless $(u\),v)\in\DlkotZ\]$.
\endpro
\Th{sl3dexp}
{\rm (Discrete exponential Selberg integral in the $\gsl_3$ case)}
$$
\alignat2
& \kern-3.14em\sum_{(u\),v)\)\in\)\DlkotZ\!\!}F(u\);v\);z_1,z_2)\,=\,
z_1^{\)k_1(k_1-1)\)\gm/2}\>z_2^{\)k_2(k_2-1)\)\gm/2}\,\x{}
\Tag{dexp3}
\\
\nn4>
{}\x\,{} &(1-z_1)^{\)(k_1-k_2)(\gm\)-\)\al\)-\)k_1\gm)}\,
(1-z_2)^{\)k_2(k_1-k_2+1)\)\gm}\,
(1-z_1z_2)^{\)k_2(\gm\)-\)\al\)-\)k_1\gm)}\,\x{}\!\} \kern-1.8em &
\\
\nn6>
&&
\Llap{{}\x\,\prod_{j=0}^{k_1-1}\,{\Gm(\al+j\gm)\,\Gm(\gm+j\)\gm)\over\Gm(\gm)}
\;\prod_{j=0}^{k_2-1}\,{\Gm(-k_1\gm+j\)\gm)\,\Gm(\gm+j\)\gm)\over\Gm(\gm)}
\kern-1.8em} &
\endalignat
$$
where the \fn/ $F$ is defined by formulae \(phi)\>--\>\(Fun)\).
\endpro
The series \(dexp3) converges for ${|\)z_1\)|<1}$, ${|\)z_2\)|<1}$ and
generic $\al\),\gm$. For ${k_2\]=0}$ formula \(dexp3) coincides with
\uugood
formula \(dexp)\).

\Section{Selberg integrals associated with $\gsl_3$}

For any nondecreasing map ${M:\lb 1\lc k_2\rb\)\to\)\lb 1\lc k_1\rb}$ \st/
$M(b)\le k_1\]-k_2\]+b$ \>for any $b=1\lc k_2$, introduce a domain $\Dkot[x,y]$
in $\Rkot\!$ with coordinates $\tko,\)\skt$ defined by the following
inequalities:
\ifMag
\vvn-.3>
$$
\gather
x\le t_{k_1}\lsym\le t_1 \le y\,,\qquad x\le s_{k_2}\lsym\le s_1 \le y\,,
\\
\nn6>
t_{M(b)}\le s_b\le t_{M(b)-1}\,,\qquad b=1\lc k_2\,.
\endgather
$$
\else
$$
x\le t_{k_1}\lsym\le t_1 \le y\,,\qquad x\le s_{k_2}\lsym\le s_1 \le y\,,
\qquad t_{M(b)}\le s_b\le t_{M(b)-1}\,,\quad b=1\lc k_2\,.
$$
\fi
Here $t_0=y$. Note that \,$\Dlkot[x,y]\)=\sum_\msm\)\Dkot[x,y]$ \>as chains.
Consider the chain
\vvn-.2>
$$
\Ckot[x,y]\,=\>\sum_\msm\)\Xkot\>\Dkot[x,y]
\vv-.6>
$$
with coefficients
\ifMag
\vvn-.2>
$$
\Xkot\>=\,\prbkt\,{\sin\)\bigl(\pi\)(k_1\]-k_2\]-M(b)+b+1\))\)\gm\bigr)
\over\sin\)\bigl(\pi\)(k_1\]-k_2\]+b\))\)\gm\bigr)}\;.
$$
\else\
$\dsize\Xkot\>=\,\prbkt\,{\sin\)\bigl(\pi\)(k_1\]-k_2\]-M(b)+b+1\))\)\gm\bigr)
\over\sin\)\bigl(\pi\)(k_1\]-k_2\]+b\))\)\gm\bigr)}\>$. \,
\vvn.2>
\fi
\vsk.2>
Introduce the \wt/ \fn/
\vvn-.8>
$$
\align
g\)(\tko;\)\skt)\, &{}=\,\Syml{\tko}\,\Syml{\skt}
\Bigl(\;\prbkt\,{1\over s_b\]-t_{b\)+k_1\]-k_2}\]}\,\Bigr)\)={}\kern-1em
\Tag{wtg}
\\
\nn4>
&{}=\,\Syml{\tko}\,\Syml{\skt}
\Bigl(\;\prbkt\,{1\over s_b\]-t_b\]}\,\Bigr)\,.
\endalign
$$
The \wt/ \fn/ is \sym/ in $\tko$ and in $\skt$ separately, and has at most
simple poles located at the hyperplanes $t_a=\)s_b$. Formula \(wtg) is
the specialization of formula \(wtf) at $\gm=0$.
\Th{sl3exp}
{\rm (\)Exponential Selberg integral in the $\gsl_3$ case)}
\vvn-.2>
$$
\alignat2
& \kern-1.76em
\int_{\Ckot[\)0,\)+\8]}\prako\>e^{-\)\bt_1t_a}\)t_a^{\>\al-1}\>
\prbkt\>e^{-\)\bt_2s_b}\,\prako\,\prbkt\,|\)t_a\]-s_b|^{\)-\)\gm}\>\x{}
\Tag{exp3}
\\
\nn4>
& {}\x\prabko\!\}(t_a\]-t_b)^{\)2\gm}\!\!\prabkt\!\}(s_a\]-s_b)^{\)2\gm}
\>g\)(\tko;\)\skt)\,\dkts\,={} \kern-2em
\\
\nn10>
&&\Llap{{}=\,\bt_1^{\)(k_1-k_2)(\gm\)-\)\al\)-\)k_1\gm)}\,
\bt_2^{\)k_2(k_1-k_2+1)\)\gm}\,
(\)\bt_1\]+\bt_2)^{\)k_2(\gm\)-\)\al\)-\)k_1\gm)}\>\x{}} \kern-2em &
\\
\nn7>
&& \Llap{{}\x\,
\prod_{j=0}^{k_1-1}\,{\Gm(\al+j\gm)\,\Gm(\gm+j\)\gm)\over\Gm(\gm)}\;
\prod_{j=0}^{k_2-1}\,{\Gm(-k_1\gm+j\)\gm)\,\Gm(\gm+j\)\gm)\over\Gm(\gm)}
\,\)\Rlap{\).}} \kern-2em &
\endalignat
$$
\endpro
\Th{sl3selb}
{\rm (\)Selberg integral in the $\gsl_3$ case)}
$$
\align
& \kern-1.24em
\int_{\Ckot[\)0,1]}\prako\,t_a^{\>\al-1}\)(1-t_a)^{\)\bt_1-1}\>
\prbkt\,(1-s_b)^{\)\bt_2-1}\,\prako\,\prbkt\,|\)t_a\]-s_b|^{\)-\)\gm}\>\x{}
\Tagg{selb3}
\\
\nn4>
& {}\x\prabko\!\}(t_a\]-t_b)^{\)2\gm}\!\!\prabkt\!\}(s_a\]-s_b)^{\)2\gm}
\>g\)(\tko;\)\skt)\,\dkts\,={}\kerm-1:0>
\\
\nn10>
& \}{}=\,\prod_{j=0}^{k_1-1}\,{\Gm(\al+j\gm)\,\Gm(\gm+j\)\gm)\over\Gm(\gm)}
\ \prod_{j=0}^{k_1-\)k_2-1}\!\!{\Gm(\)\bt_1\]+j\)\gm)\over
\Gm\bigl(\al+\bt_1\]+(2\)k_1\]-k_2\]-2-j)\)\gm\bigr)}\;\x{}
\\
\nn7>
& {}\x\,\prod_{j=0}^{k_2-1}\,{\Gm(\bt_2\]+j\)\gm)\,
\Gm(\bt_1\]+\bt_2\]-1-\gm+j\gm)\,\Gm(-k_1\gm+j\)\gm)\,\Gm(\gm+j\)\gm)\over
\Gm\bigl(\bt_2\]+(2\)k_2\]-k_1\]-2-j\))\gm\bigr)\,
\Gm\bigl(\al+\bt_1\]+\bt_2\]-1+(\)k_1\]+k_2\]-3-j\))\gm\bigr)\,
\Gm(\gm)}\;.\kerm-3:-2>
\\
\endalign
$$
\endpro
\Th{sl3selb0}
{\rm (\)Another Selberg integral in the $\gsl_3$ case)}
$$
\align
& \kern-1.24em
\int_{\Ckot[\)0,1]}\prako\,t_a^{\>\al-1}\)(1-t_a)^{\)\bt_1-1}\>
\prbkt\,(1-s_b)^{\)\bt_2-1}\,\prako\,\prbkt\,|\)t_a\]-s_b|^{\)-\)\gm}\>\x{}
\Tagg{selb30}
\\
\nn4>
& \hp{\kern-3.27em
\int_{\Ckot[\)0,1]}\prako\,t_a^{\>\al-1}\)(1-t_a)^{\)\bt_1-1}}\x
\prabko\!\}(t_a\]-t_b)^{\)2\gm}\!\!\prabkt\!\}(s_a\]-s_b)^{\)2\gm}\,\dkts\,={}\kerm-1:0>
\\
\nn10>
& \}{}=\,\prod_{j=0}^{k_1-1}\,{\Gm(\al+j\gm)\,\Gm(\gm+j\)\gm)\over\Gm(\gm)}
\ \prod_{j=0}^{k_1-\)k_2-1}\!\!{\Gm(\)\bt_1\]+j\)\gm)\over
\Gm\bigl(\al+\bt_1\]+(2\)k_1\]-k_2\]-2-j)\)\gm\bigr)}\;\x{}
\\
\nn7>
& {}\x\,\prod_{j=0}^{k_2-1}\,{\Gm(\bt_2\]+j\)\gm)\,
\Gm(\bt_1\]+\bt_2\]-\gm+j\gm)\,\Gm(1-k_1\gm+j\)\gm)\,\Gm(\gm+j\)\gm)\over
\Gm\bigl(\bt_2\]+1+(2\)k_2\]-k_1\]-2-j\))\gm\bigr)\,
\Gm\bigl(\al+\bt_1\]+\bt_2\]+(\)k_1\]+k_2\]-3-j\))\gm\bigr)\,
\Gm(\gm)}\;.\kerm-3:-2>
\endalign
$$
\endpro
\nt
Integrals \(exp3)\>--\>\(selb30) converge if \>${\Re\al>0}$, ${\Re\bt_1>0}$,
${\Re\bt_2>0}$, ${\Re\gm<0}$ and $|\]\Re\gm\>|$ is sufficiently small.
For $k_2\]=0$ formulae \(exp3) and \(selb3)\), \(selb30) coincide \resp/
with formulae \(exp) and \(selb)\).
\vsk.2>
There is another description of the integration chains $\Ckot\}$ in
Theorems~\[sl3exp]\>--\,\[sl3selb0]\). We explain it for the exponential
Selberg integral \(exp3)\). For the Selberg integrals \(selb3) and \(selb30)
the description is similar.
\vsk.1>
Consider a collection of simple curves $\Cc_1\lc\Cc_{k_1\]+\)k_2}$ going
around the real positive semiline $\Rp$ like the curve in the picture:
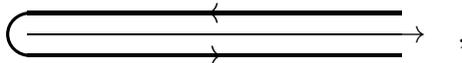
\vvn.4>
$$
\roop\Rlap{\quad,}
\vv.6>
$$
and \st/ $\Cc_i$ is between $\Rp$ and $\Cc_j$ for any $i<j$.
Let $G(\tko;\)\skt)$ be the integrand in \lhs/ of formula \(exp3) and let
$$
\align
\Gt(\tko;{}& \)\skt)=\prako\,e^{-\)\bt_1t_a}\)t_a^{\>\al-1}\>
\prbkt e^{-\)\bt_2s_b}\,\prako\,\prbkt\,(s_b\]-t_a)^{\)-\)\gm}\>\x{}
\\
\nn5>
& \!\]{}\x\prabko\!\}(t_a\]-t_b)^{\)2\gm}\!\!\prabkt\!\}(s_a\]-s_b)^{\)2\gm}
\>g\)(\tko;\)\skt)
\endalign
$$
(no absolute value of ${s_b\]-t_a}$)\). Fix a branch of the \fn/
$\Gt(\tko;\)\skt)$ on ${\Cc_1\}\lx\)\Cc_{k_1+\)k_2}}$ in the following way:
at the point, where all the \var/s $\tko,\)\skt$ are real
negative, one has \>$\arg\)t_a=\arg\)({s_b\]-t_a})=\)\pi$ for any $a\),b$,
and \>$\arg\)({t_a\]-t_b})=\arg\)({s_a\]-s_b})=\)0$ for ${a<b}$.
Equivalently, if ${\Im t_a>0}$, ${\Im s_b>0}$ for any $a\),b$, and
$0<\)\Re t_{k_1}\}\lsym<\)\Re t_1<\)\Re s_{k_2}\}\lsym<\)\Re s_1$, then
$\arg\)t_a$, ${\arg\)(s_b\]-t_a)}$ for any $a\),b$, and ${\arg\)(t_a\]-t_b)}$,
${\arg\)(s_a\]-s_b)}$ for ${a<b}$ are between $-\)\pi/2$ and $\pi/2$.
\Th{chains}
$$
\align
&\kern-1.76em\int_{\Cc_1\]\lx\)\Cc_{k_1+\)k_2}\!\!}
\!\Gt(\tko;\)\skt)\,\dkts\,={}
\\
& {}=\,N_{k_1}(\al\),\gm)\,N_{k_2}(-k_1\gm\),\gm)
\int_{\Ckot[\)0,\)+\8]}\!\!\!G(\tko;\)\skt)\,\dkts\,.\kern-1em
\\
\cnn.1>
\endalign
$$
where\hfill
$\dsize N_k(\al\),\gm)\,=\,
\prod_{j=0}^{k-1}\;{2\)i\>e^{\)\pi i\al}\sin\)\bigl(\pi\)(\al+j\)\gm)\bigr)
\)\sin\)\bigl(\pi(\gm+j\)\gm)\bigr)\over\sin\)(\pi\gm)}\;.$\hfill\hp{w}
\endpro
\vsk.4>
\nt The proof is straightforward.

\Section{On proofs of Theorems \[sl3dexp]\), \[sl3exp]\,--\,\)\[sl3selb0]}

All four types of new Selberg integrals \(dexp3)\), \(exp3)\>--\>\(selb30) are
related to \rep/ theory of the Lie algebra $\gsl_3$. They appear as coordinate
\fn/s of \hgeom/ \sol/s of the \KZ/ and dynamical \eq/s associated with
$\gsl_3$.
\vsk.1>
Let $\al_1\),\al_2$ be simple roots of $\gsl_3$, and let $\om_1\),\om_2$ be
the fundamental \wt/s. Denote by $L_\la$ the \irr/ \$\gsl_3\)$-module of \hw/
$\la$. Series \(dexp3) appears as the coordinate \fn/ of the \hgeom/ \sol/ of
the dynamical \difl/ \eq/s from \cite{TV3} with values in the \wt/ subspace
$L_\la\)[\)\la-k_1\al_1\]-k_2\)\al_2\)]$ for ${\la\in\C\>\om_1}$, and integral
\(exp3) appears as the coordinate \fn/ of the \hgeom/ \sol/ of the dynamical
\difl/ \eq/s from \cite{FMTV} with values in the same \wt/ subspace,
see~\cite{FMTV, formula (8)\)}\). Integrals \(selb3)\), \(selb30) appear
as coordinate \fn/s of the \hgeom/ \sol/ of the \KZ/ \difl/ \eq/s
with values in the space
$\Sing(L_{\la_1}\!\ox L_{\la_2})\)[\)\la_1\]+\la_2-k_1\al_1\]-k_2\)\al_2\)]$
of singular \wt/ vectors for $\la_1\}\in\C\>\om_1$, see \cite{SV}\),
and Chapter~12 in \cite{V}\).
\vsk.1>
To prove Theorem~\[sl3dexp] one considers expressions in both sides of formula
\(dexp3) as \fn/s of the \var/s $z_1\), z_2$. The dynamical \difl/ \eq/s in
the case in question take the form
\vvn-.2>
$$
\gather
{\der\Psi\over\der z_1}\;=\,\biggl(\>{k_1(k_1-1)\)\gm\over 2\)z_1}\,+\,
{(k_1-k_2)\>(\al-\gm+k_1\gm)\over 1-z_1}\,+\,
{z_2\)k_2(\al-\gm+k_1\gm)\over 1-z_1z_2}\>\biggr)\,\Psi\,,
\\
\nn7>
{\der\Psi\over\der z_2}\;=\,\biggl(\>{k_2(k_2-1)\)\gm\over 2\)z_1}\,-\,
{k_2(k_1-k_2+1)\)\gm\over 1-z_2}\,+\,{z_1\)k_2(\al-\gm+k_1\gm)\over 1-z_1z_2}
\>\biggr)\,\Psi\,.
\endgather
$$
The series in \lhs/ solves this system and, hence, equals the product in \rhs/
up to a factor which does not depend on $z_1\), z_2$. Since the values of both
expressions coincide at $z_1=\)z_2\)=0$, the proportionality factor equals $1$,
which completes the proof.
\vsk.1>
Theorem~\[sl3exp] is derived from Theorem~\[sl3dexp] by taking the limit
$z_1\]\to 1$, $z_2\]\to 1$. Namely, set $z_1=e^{\)-\eps\bt_1}\}$,
$z_2=e^{\)-\eps\bt_2}\}$, $u_a=t_a/\eps$, $v_b=s_b/\eps$, and take the limit
$\eps\to 0$. If \>${\Re\al>0}$, ${\Re\bt_1>0}$, ${\Re\bt_2>0}$, ${\Re\gm<0}$
and $|\]\Re\gm\>|$ is sufficiently small, then using a corollary of
the Stirling formula:
\vvn-.3>
$$
{\Gm(x+c)\over\Gm(x+d)}\;=\,x^{\)c\)-\)d}\)\ono\,,\qquad \Re x\to+\)\8\,,
$$
one shows that formula \(dexp3) tends to formula \(exp3) in the limit.
This is quite similar to getting formula \(exp) from formula \(dexp)
in the limit $z\to 1$.
\vsk.1>
To obtain Theorem~\[sl3selb] one shows that expressions in both sides of
formula \(dexp3) obey the same system of \deq/s \wrt/ the shifts of the \var/s
\)$\bt_1\]\to\bt_1\]+1$ and \)$\bt_2\]\to\bt_2\]+1$, see \cite{MaV}\); those
\deq/s are the dynamical \deq/s from \cite{TV2}\). Therefore, those expressions
are proportional up to a periodic \fn/ of $\bt_1\),\bt_2$. The periodic \fn/
can be found by comparing asymptotics as $\Re\bt_1\to\)+\)\8$,
$\Re\bt_2\to\)+\)\8$. In this limit formula \(selb3) reduces to formula
\(exp3), which shows that the periodic \fn/ is the constant equal to $1$.
\vsk.1>
The proof of Theorem~\[sl3selb0] is based on formula \(selb3)\).
It is outlined in the next section, \cf. \(J0k)\).

\Section{Further generalizations of Selberg type integrals}

Denote by $S(\al\),\bt\))$ the Selberg integral \(selb)\).
Consider the following integrals:
\vvn-.2>
$$
I_{\)l}\,=\!\int_{\Dlk[\)0,1]}
\prak\,t_a^{\>\al-1}\)(1-t_a)^{\)\bt-1}\!\!\prabk\!\}(t_a\]-t_b)^{\)2\gm}\;
\Syml{\tk}\!\]\Bigl(\,\)\tprod_{a=1}^l t_a\tprod_{b=l+1}^k(1-t_a)\Bigr)\,
\dkt\,,
\vv-.2>
$$
$l=0\lc k$. In particular,
\,$I_0=\)S(\al\),\bt+1)$ \>and \,$I_k=\)S(\al+1,\bt\))$.
\Prop{more}
For any \,$l=0\lc k$ \>one has
\vvn-.2>
$$
I_{\)l}\,=\,\prod_{i=0}^{l-1}\,{\al+(k-1-i\))\)\gm\over\bt+i\)\gm}\,\,
\prod_{j=0}^{k-1}\,{\Gm(\al+j\)\gm)\,\Gm(1+\bt+j\)\gm)\,\Gm(\gm+j\gm)
\over\Gm\bigl(1+\al+\bt+(2\)k-2-j\))\gm\bigr)\,\Gm(\gm)}\;.
$$
\endpro
\nt
The proof is straightforward using Aomoto's formula \cite{A1}\):
\vvn-.2>
$$
\align
\int_{\Dlk[\)0,1]}{} &
\prak\,t_a^{\>\al-1}\)(1-t_a)^{\)\bt-1}\!\!\prabk\!\}(t_a\]-t_b)^{\)2\gm}\;
\Syml{\tk}\!\](t_1\ldots t_l)\,\dkt\,={}
\\
\nn2>
{}=\,{} & \prod_{i=0}^{\Cph{a=1}{l-1}}
\,{\al+(k-1-i\))\)\gm\over\al+\bt+(2\)k-2-i\))\gm}\,\,
\prod_{j=0}^{k-1}\,{\Gm(\al+j\)\gm)\,\Gm(\)\bt+j\)\gm)\,\Gm(\gm+j\gm)
\over\Gm\bigl(\al+\bt+(2\)k-2-j\))\gm\bigr)\,\Gm(\gm)}\;.
\endalign
$$
Proposition \[more] is \eqv/ to the formula for the Selberg integral
and linear relations for the integrals $I_0\lc I_k$:
\vvn.2>
$$
\bigl(\al+(k-l-1)\)\gm\bigr)\>I_{\)l}\>=\,(\bt+l\gm)\>I_{\)l+1}\,,\qqq
l=0\lc k-1\,.\kern-3em
\vv.3>
$$
These relations correspond to the fact that the vector
${\sum_{a=0}^k (-1)^a\>I_a\>f^{\)k-a}v_1\ox f^a v_2}$ is a singular vector in
the tensor product ${L_{\ell_1}\!\ox L_{\ell_2}}\}$ of \hw/ \$\gsl_2\)$-modules
with \hw/s \>$\ell_1=-\)\al/\gm$, $\ell_2=-\)\bt/\gm$, see~\cite{V}\).
Recall that a vector $v$ is called singular if $e\)v=0$. Here $e\>,\]f$
are the standard generators of the opposite nilpotent subalgebras of the Lie
algebra $\gsl_2$, and $v_1\),\)v_2$ are the \hwv/s (standard cyclic vectors)
of the respective \$\gsl_2\)$-modules $L_{\ell_1},\)L_{\ell_2}\}$.
\vsk.2>
Fix nonnegative integers $k_1\),k_2$ \st/ $k_1\ge k_2$.
Denote by $R\)(\al\),\)\bt_1,\)\bt_2)$ the integral~\(selb3)\).
Consider the master \fn/
\ifMag
$$
\alignat2
& \Om\)(\tko;\)\skt;\al\),\)\bt_1,\)\bt_2)\,={}
\\
\nn7>
&\!\}{}=\,\prako\,t_a^{\>\al-1}\)(1-t_a)^{\)\bt_1-1}\prbkt(1-s_b)^{\)\bt_2-1}
\,\prako\,\prbkt\,|\)t_a\]-s_b|^{\)-\)\gm}\)\x{}\!\!&&
\\
\nn6>
&&\Llap{
{}\x\prabko\!\}(t_a\]-t_b)^{\)2\gm}\!\!\prabkt\!\}(s_a\]-s_b)^{\)2\gm}}\,.
\\
\cnn.2>
\endalignat
$$
\else
$$
\align
& \Om\)(\tko;\)\skt;\al\),\)\bt_1,\)\bt_2)\,={}
\\
\nn5>
& \!\}{}=\,\prako\,t_a^{\>\al-1}\)(1-t_a)^{\)\bt_1-1}
\prbkt (1-s_b)^{\)\bt_2-1}\,\prako\,\prbkt\,|\)t_a\]-s_b|^{\)-\)\gm}
\prabko\!\}(t_a\]-t_b)^{\)2\gm}\!\!\prabkt\!\}(s_a\]-s_b)^{\)2\gm}\,.
\endalign
$$
\fi
Say that a triple of nonnegative integers $(\)l_1\),\)l_2\),\)m)$ is
admissible if ${l_1\le k_1\]-k_2\]+l_2}$, \>${l_2\le k_2}$ \>and
\>${m\le\min\)(\)\l_1\),l_2\))}$. For any admissible triple
$l_1\),\)l_2\),\)m$ introduce the \fn/s
$$
\align
h_{\)l_1\],\>l_2,\)m} &{}(\tko;\)\skt)\,={}
\\
\nn3>
\;{}=\>{} & \)\Syml{\tko}\,\Syml{\skt}\Bigl(\;
\prod_{a=1}^{l_1}\,t_a\prod_{a=l_1+1}^{k_1}(1-t_a)\>
\prod_{b=1}^{m}\,{1-s_b\over s_b\]-t_b\]}\;
\prod_{b=l_2+1}^{k_2}\,{1-s_b\over s_b\]-t_{b\)+k_1\]-k_2}\]}\,\Bigr)\,,
\\
\nn9>
\hti_{\)l_1\],\>l_2,\)m} &{}(\tko;\)\skt)\,={}
\\
\nn3>
\;{}=\>{} & \)\Syml{\tko}\,\Syml{\skt}\Bigl(\;
\prod_{a=1}^{l_1}\,t_a\prod_{a=l_1+1}^{k_1}(1-t_a)\;
\prod_{b=1}^{m}\,{1-t_b\over s_b\]-t_b\]}\;
\prod_{b=l_2+1}^{k_2}\,{1-s_b\over s_b\]-t_{b\)+k_1\]-k_2}\]}\,\Bigr)\,,
\endalign
$$
and the integrals
\ifMag\else\vvn-.7>\fi
$$
\align
& J_{\)l_1\],\>l_2,\)m}\,=\!\int_{\Ckot[\)0,1]}\kern-.55em
\Om(t\);s\);\al\),\)\bt_1,\)\bt_2)\,\>h_{\)l_1\],\>l_2,\)m}(t\);s\))\,\dkts\,,
\\
\nn7>
& \Jt_{\)l_1\],\>l_2,\)m}\,=\!\int_{\Ckot[\)0,1]}\kern-.55em
\Om(t\);s\);\al\),\)\bt_1,\)\bt_2)\,\>\hti_{\)l_1\],\>l_2,\)m}(t\);s\))\,
\dkts\,.
\\
\cnn-.8>
\nngood
\cnn.8>
\endalign
$$
In particular, \,$J_{\)0,\)0,\)0}\)=\)
\Jt_{\)0,\)0,\)0}\)=\)R\)(\al\),\)\bt_1\]+1,\)\bt_2\]+1)$.
One can also observe that
\vvn.3>
$$
J_{\)0,\>l,\)0}(\al+1,\)\bt_1,\)\bt_2)\,=\,
J_{k_1\],\)k_2,\)k_2-\)l}(\al\),\)\bt_1\]+1,\)\bt_2)\,.
\vv.2>
\Tag{JJl}
$$
Here we indicate the dependence on the parameters $\al\),\)\bt_1\),\bt_2$
explicitly. In addition, the integral
$$
J_{\)0,\)k_2,\)0}\>=\>\Jt_{\)0,\)k_2,\)0}\,=\!
\int_{\Ckot[\)0,1]}\kern-.55em\Om(t\);s\);\al,\)\bt_1+1,\)\bt_2)\,\dkts
\vv-.2>
\Tag{JJ0}
$$
is expressed in terms of the master \fn/ only with no extra \raf/s involved.
\Th{JJJ}
Let $(\)l_1\),\)l_2\),\)m)$ be an admissible triple.
If \,${l_1<k_1\]-k_2\]+l_2}$, \>then
\ifMag
$$
\align
& \bigl(\al+(k_1\]-l_1\]-1)\)\gm\bigr)\>J_{\)l_1\],\>l_2,\)m}\>={}
\\
\nn4>
& \>{}=\,\bigl(\)\bt_1\]+(\)l_1\]-l_2\]+m)\)\gm\bigr)\>J_{\)l_1+1,\>l_2,\)m}
\)-\>(l_2\]-m)\)\gm\>J_{\)l_1+1,\>l_2,\)m+1}\,,
\endalign
$$
\else
$$
\bigl(\al+(k_1\]-l_1\]-1)\)\gm\bigr)\>J_{\)l_1\],\>l_2,\)m}\>=\,
\bigl(\)\bt_1\]+(\)l_1\]-l_2\]+m)\)\gm\bigr)\>J_{\)l_1+1,\>l_2,\)m}
\)-\>(l_2\]-m)\)\gm\>J_{\)l_1+1,\>l_2,\)m+1}\,,
$$
\fi
and if \,$m<l_2<k_2$, \>then
\ifMag
$$
\align
& (k_2\]-k_1\]+l_1\]-l_2)\)\gm\>J_{\)l_1\],\>l_2-1,\)m}\>={}
\Tag{J1}
\\
\nn4>
& \){}=\,\bigl(\)\bt_2\]+(\)l_2\]-m-1)\)\gm\bigr)\>J_{\)l_1\],\>l_2,\)m}
\)+\>(l_1\]-m)\)\gm\>J_{\)l_1\],\>l_2,\)m+1}\,.
\endalign
$$
\else
$$
(k_2\]-k_1\]+l_1\]-l_2)\)\gm\>J_{\)l_1\],\>l_2-1,\)m}\>=\,
\bigl(\)\bt_2\]+(\)l_2\]-m-1)\)\gm\bigr)\>J_{\)l_1\],\>l_2,\)m}
\)+\>(l_1\]-m)\)\gm\>J_{\)l_1\],\>l_2,\)m+1}\,.
\Tag{J1}
$$
\fi
Similarly, if \,${l_1<k_1\]-k_2\]+l_2}$, \>then
\ifMag
$$
\align
& \bigl(\al+(k_1\]-l_1\]-1)\)\gm\bigr)\>\Jt_{\)l_1\],\>l_2,\)m}\>={}
\\
\nn4>
& \>{}=\,
\bigl(\)\bt_1\]+(\)l_1\]-m)\)\gm\bigr)\>\Jt_{\)l_1+1,\>l_2,\)m}
\)-\>(l_2\]-m)\)\gm\>\Jt_{\)l_1+1,\>l_2,\)m+1}\,,
\endalign
$$
\else
$$
\bigl(\al+(k_1\]-l_1\]-1)\)\gm\bigr)\>\Jt_{\)l_1\],\>l_2,\)m}\>=\,
\bigl(\)\bt_1\]+(\)l_1\]-m)\)\gm\bigr)\>\Jt_{\)l_1+1,\>l_2,\)m}
\)-\>(l_2\]-m)\)\gm\>\Jt_{\)l_1+1,\>l_2,\)m+1}\,,
$$
\fi
and if \,$m<l_2<k_2$, \>then
\ifMag
$$
\align
& (k_2\]-k_1\]+l_1\]-l_2)\)\gm\>\Jt_{\)l_1\],\>l_2-1,\)m}\>={}
\Tag{Jt1}
\\
\nn4>
&\){}=\,\bigl(\)\bt_2\]+(\)l_2\]-l_1\]+m-1)\)\gm\bigr)\>\Jt_{\)l_1\],\>l_2,\)m}
\)+\>(l_1\]-m)\)\gm\>\Jt_{\)l_1\],\>l_2,\)m+1}\,.
\endalign
$$
\else
$$
(k_2\]-k_1\]+l_1\]-l_2)\)\gm\>\Jt_{\)l_1\],\>l_2-1,\)m}\>=\,
\bigl(\)\bt_2\]+(\)l_2\]-l_1\]+m-1)\)\gm\bigr)\>\Jt_{\)l_1\],\>l_2,\)m}
\)+\>(l_1\]-m)\)\gm\>\Jt_{\)l_1\],\>l_2,\)m+1}\,.
\Tag{Jt1}
$$
\fi
\endpro
\ifMag\ifUS\else\nt\fi\fi
In the \eq/s of Theorem~\[JJJ] some of the triples
$(\){l_1+1}\),\)l_2\),\){m+1}\))$, \>$(\)l_1\),\){l_2-1}\),\)m)$,
\>$(\)l_1\),\)l_2\),\){m+1}\))$ may be inadmissible. In that case
the corresponding coefficient:
\ifMag
\vvn-.1>
$$
l_2\]-m\,,\quad k_2\]-k_1\]+l_1\]-l_2\,,\quad l_1\]-m\,,
\vv-.1>
$$
\else
${l_2\]-m}$, \>${k_2\]-k_1\]+l_1\]-l_2}$, \>${l_1\]-m}$,
\fi
\resp/, equals zero, and the integral with an inadmissible index does not enter
the equation.
\vsk.5>
Let $f_1\),f_2\in\gsl_3$ be root vectors corresponding to simple roots, and
let $\om_1\),\)\om_2$ be the fundamental \wt/s. It is known that the vectors
$$
\sum_{l_1\],\>l_2,\)m}\](-1)^{l_1}\)J_{\)l_1\],\>l_2,\)m}\,
{f_1^{\)k_1\]-k_2\]-l_1\]+\)l_2}\>[\)f_1\),\]f_2\)]^{\)k_2\]-l_2}\>v_1\ox
f_1^{\>l_1\]-m}\>[\)f_1\),\]f_2\)]^{\)m}\)f_2^{\>l_2\]-m}\>v_2\over
(k_1\]-k_2\]-l_1\]+l_2)!\,(k_2\]-l_2)!\,(l_1\]-m)!\,\)m\)!\,(l_2\]-m)!}
\vv-.1>
$$
and
\vvn-.5>
$$
\sum_{l_1\],\>l_2,\)m}\](-1)^{l_1}\)\Jt_{\)l_1\],\>l_2,\)m}\,
{f_1^{\)k_1\]-k_2\]-l_1\]+\)l_2}\>[\)f_1\),\]f_2\)]^{\)k_2\]-l_2}\>v_1\ox
f_2^{\>l_2\]-m}\>[\)f_1\),\]f_2\)]^{\)m}\)f_1^{\>l_1\]-m}\>v_2\over
(k_1\]-k_2\]-l_1\]+l_2)!\,(k_2\]-l_2)!\,(l_1\]-m)!\,\)m\)!\,(l_2\]-m)!}
$$
are singular vectors, that is, they are annihilated by the action of
the opposite root vectors $e_1\), e_2\in\gsl_3$, in the tensor product
${L_{\la_1}\!\ox L_{\la_2}}\}$ of \$\gsl_3$-modules with \hw/s
\ifMag
\ifUS\else\vvn-.1>\fi
$$
\la_1\>=\,-\)\om_1\al/\gm\,,\qquad
\la_2\>=\,-\)\om_1\bt_1/\gm\)-\om_2\bt_2/\gm
\ifUS\else\vv.2>\fi
$$
\else
$\la_1=-\)\om_1\al/\gm$, $\la_2=-\)\om_1\bt_1/\gm\)-\om_2\bt_2/\gm$
\fi
and \hw/ vectors $v_1\),\)v_2$, see~\cite{V}\).
This fact implies the relations of Theorem~\[JJJ].
\ifUS\else\vskmgood-.7:.7> \fi
\vsk.2>
For generic values of $\al\),\)\bt_1\),\bt_2\),\gm$
\)one can solve \eq/s of Theorem~\[JJJ] and express all the integrals
$J_{\)l_1\],\>l_2,\)m}$ and $\Jt_{\)l_1\],\>l_2,\)m}$ via the known integral
\ifMag
$$
J_{\)0,\)0,\)0}\,=\,\Jt_{\)0,\)0,\)0}\,=\,R\)(\al\),\)\bt_1\]+1,\)\bt_2\]+1)\,.
$$
\else
$J_{\)0,\)0,\)0}\)=\)\Jt_{\)0,\)0,\)0}\)=\)R\)(\al\),\)\bt_1\]+1,\)\bt_2\]+1)$.
\fi
For instance, the integrals $J_{\)0,\>l,\)0}=\Jt_{\)0,\>l,\)0}$ can be easily
computed using relations \(J1)\):
$$
J_{\)0,\>l,\)0}\,=\,(-1)^{\)l}\,\,\prod_{i=0}^{l-1}\;
{(k_1\]-k_2\]+1+i\))\)\gm\over\bt_2\]+i\)\gm}\ \,J_{\)0,\)0,\)0}\,.
\ifMag\ifUS\vv-.5>\else\vv-.8>\fi\fi
$$
In particular,
\ifMag
\ifUS
\vvn-.4>
$$
\alignat2
& J_{\)0,\)k_2,\)0}\,=\,
\Tagg{J0k}
\\
\nn10>
& \quad{}=\prod_{j=0}^{k_1-1}\,{\Gm(\al+j\gm)\,\Gm(\gm+j\)\gm)\over\Gm(\gm)}
\ \prod_{j=0}^{k_1-\)k_2-1}\!\!{\Gm(\)\bt_1\]+1+j\)\gm)\over
\Gm\bigl(\al+\bt_1\]+1+(2\)k_1\]-k_2\]-2-j)\)\gm\bigr)}\;\x{}\;\kern-3.4em
\\
\nn4>
&& \Llap{{}\x\,\prod_{j=0}^{k_2-1}\,{\Gm(\bt_2\]+j\)\gm)\,
\Gm(\bt_1\]+\bt_2\]+1-\gm+j\gm)\,\Gm(1-k_1\gm+j\)\gm)\,\Gm(\gm+j\)\gm)\over
\Gm\bigl(\bt_2\]+1+(2\)k_2\]-k_1\]-2-j\))\gm\bigr)\,
\Gm\bigl(\al+\bt_1\]+\bt_2\]+1+(\)k_1\]+k_2\]-3-j\))\gm\bigr)\,
\Gm(\gm)}\;.\kern-3.4em}
\endalignat
$$
\else
$$
\gather
{\align
& J_{\)0,\)k_2,\)0}\,=\,
\Tagg{J0k}
\\
\nn10>
& \!{}=\prod_{j=0}^{k_1-1}\,{\Gm(\al+j\gm)\,\Gm(\gm+j\)\gm)\over\Gm(\gm)}
\ \prod_{j=0}^{k_1-\)k_2-1}\!\!{\Gm(\)\bt_1\]+1+j\)\gm)\over
\Gm\bigl(\al+\bt_1\]+1+(2\)k_1\]-k_2\]-2-j)\)\gm\bigr)}\;\x{}\;\kern-3em
\endalign}
\\
\nn4>
\Rline{{}\x\,\prod_{j=0}^{k_2-1}\,{\Gm(\bt_2\]+j\)\gm)\,
\Gm(\bt_1\]+\bt_2\]+1-\gm+j\gm)\,\Gm(1-k_1\gm+j\)\gm)\,\Gm(\gm+j\)\gm)\over
\Gm\bigl(\bt_2\]+1+(2\)k_2\]-k_1\]-2-j\))\gm\bigr)\,
\Gm\bigl(\al+\bt_1\]+\bt_2\]+1+(\)k_1\]+k_2\]-3-j\))\gm\bigr)\,
\Gm(\gm)}\;.}
\endgather
$$
\fi
\else
\vvn-.3>
$$
\align
J_{\)0,\)k_2,\)0}\, &{}=\,
\prod_{j=0}^{k_1-1}\,{\Gm(\al+j\gm)\,\Gm(\gm+j\)\gm)\over\Gm(\gm)}
\ \prod_{j=0}^{k_1-\)k_2-1}\!\!{\Gm(\)\bt_1\]+1+j\)\gm)\over
\Gm\bigl(\al+\bt_1\]+1+(2\)k_1\]-k_2\]-2-j)\)\gm\bigr)}\;\x{}
\Tagg{J0k}
\\
\nn4>
& \>{}\x\,\prod_{j=0}^{k_2-1}\,{\Gm(\bt_2\]+j\)\gm)\,
\Gm(\bt_1\]+\bt_2\]+1-\gm+j\gm)\,\Gm(1-k_1\gm+j\)\gm)\,\Gm(\gm+j\)\gm)\over
\Gm\bigl(\bt_2\]+1+(2\)k_2\]-k_1\]-2-j\))\gm\bigr)\,
\Gm\bigl(\al+\bt_1\]+\bt_2\]+1+(\)k_1\]+k_2\]-3-j\))\gm\bigr)\,
\Gm(\gm)}\;.\kern-2.6em
\endalign
$$
\fi
After replacement $\bt_1\to\bt_1\]-1$ the last formula yields formula
\(selb30)\), \cf. \(JJ0)\). Similarly, the integrals $\Jt_{\)k_1\],\)k_2,\)m}$
can be computed using relations \(Jt1) and the integral
$J_{k_1\],\)k_2,\)0}\)=\)\Jt_{k_1\],\)k_2,\)0}$:
\ifMag
\ifUS
\vvn-.8>
$$
\gather
\Jt_{k_1\],\)k_2,\)m}\,=\,(-1)^{\)m}\,\)\prod_{i=0}^{m-1}\,
{\bt_2\]+(k_2\]-k_1\]-1+i\))\)\gm\over(k_1\]-i\))\)\gm}
\ \,J_{k_1\],\)k_2,\)0}\,,
\\
\nn10>
{\align
J_{k_1\],\)k_2,\)0}\,=\,
\prod_{j=0}^{k_1-1}\,{\Gm(\al+1+j\gm)\,\Gm(\gm+j\)\gm)\over\Gm(\gm)}
\ \prod_{j=0}^{k_1-\)k_2-1}\!\!{\Gm(\)\bt_1\]+j\)\gm)\over
\Gm\bigl(\al+\bt_1\]+1+(2\)k_1\]-k_2\]-2-j)\)\gm\bigr)}\;\x{} &
\\
\nn4>
{}\x\,\prod_{j=0}^{k_2-1}\,{\Gm(\bt_2\]+j\)\gm)\,
\Gm(\bt_1\]+\bt_2\]-\gm+j\gm)\,\Gm(1-k_1\gm+j\)\gm)\,\Gm(\gm+j\)\gm)\over
\Gm\bigl(\bt_2\]+1+(2\)k_2\]-k_1\]-2-j\))\gm\bigr)\,
\Gm\bigl(\al+\bt_1\]+\bt_2\]+1+(\)k_1\]+k_2\]-3-j\))\gm\bigr)\,\Gm(\gm)}\;&,
\endalign}
\endgather
$$
\else
\vvn->
$$
\gather
\Jt_{k_1\],\)k_2,\)m}\,=\,(-1)^{\)m}\,\)\prod_{i=0}^{m-1}\,
{\bt_2\]+(k_2\]-k_1\]-1+i\))\)\gm\over(k_1\]-i\))\)\gm}
\ \,J_{k_1\],\)k_2,\)0}\,,
\\
\nn10>
\Lline{J_{k_1\],\)k_2,\)0}\,=\,
\prod_{j=0}^{k_1-1}\,{\Gm(\al+1+j\gm)\,\Gm(\gm+j\)\gm)\over\Gm(\gm)}
\ \prod_{j=0}^{k_1-\)k_2-1}\!\!{\Gm(\)\bt_1\]+j\)\gm)\over
\Gm\bigl(\al+\bt_1\]+1+(2\)k_1\]-k_2\]-2-j)\)\gm\bigr)}\;\x{}}
\\
\nn4>
\Rline{\>{}\x\,\prod_{j=0}^{k_2-1}\,{\Gm(\bt_2\]+j\)\gm)\,
\Gm(\bt_1\]+\bt_2\]-\gm+j\gm)\,\Gm(1-k_1\gm+j\)\gm)\,\Gm(\gm+j\)\gm)\over
\Gm\bigl(\bt_2\]+1+(2\)k_2\]-k_1\]-2-j\))\gm\bigr)\,
\Gm\bigl(\al+\bt_1\]+\bt_2\]+1+(\)k_1\]+k_2\]-3-j\))\gm\bigr)\,\Gm(\gm)}
\Rlap{\;,}}
\\
\cnn.1>
\endgather
$$
\fi
\else
\vvn-.2>
$$
\gather
\Jt_{k_1\],\)k_2,\)m}\,=\,(-1)^{\)m}\,\)\prod_{i=0}^{m-1}\,
{\bt_2\]+(k_2\]-k_1\]-1+i\))\)\gm\over(k_1\]-i\))\)\gm}
\ \,J_{k_1\],\)k_2,\)0}\,,
\\
\nn12>
{\align
J_{k_1\],\)k_2,\)0}\, &{}=\,
\prod_{j=0}^{k_1-1}\,{\Gm(\al+1+j\gm)\,\Gm(\gm+j\)\gm)\over\Gm(\gm)}
\ \prod_{j=0}^{k_1-\)k_2-1}\!\!{\Gm(\)\bt_1\]+j\)\gm)\over
\Gm\bigl(\al+\bt_1\]+1+(2\)k_1\]-k_2\]-2-j)\)\gm\bigr)}\;\x{}
\\
\nn4>
&\>{}\x\,\prod_{j=0}^{k_2-1}\,{\Gm(\bt_2\]+j\)\gm)\,
\Gm(\bt_1\]+\bt_2\]-\gm+j\gm)\,\Gm(1-k_1\gm+j\)\gm)\,\Gm(\gm+j\)\gm)\over
\Gm\bigl(\bt_2\]+1+(2\)k_2\]-k_1\]-2-j\))\gm\bigr)\,
\Gm\bigl(\al+\bt_1\]+\bt_2\]+1+(\)k_1\]+k_2\]-3-j\))\gm\bigr)\,\Gm(\gm)}\;,
\endalign}
\endgather
$$
\fi
the last formula following from \(JJl)\), \(J0k).

\section{Errata}
The published version of the paper (\>\LMP/ {\bf 65} (2003), no\point 3,
173\)\~185\)) contains a sign misprint in the formulae for the \fn/s
$h_{\)l_1\],\>l_2,\)m}\}$ and $\smash{\hti_{\)l_1\],\>l_2,\)m}}\}$ in
Section~5. This misprint is corrected in the present version.
\vsk.1>
We thank O\&Warnaar for a valuable remark.

\myRefs
\widest{EFK}

\ref\Key A1
\by \Aomoto/
\paper
Jacobi \pol/s associated with Selberg integrals
\jour \SIAM/ \yr 1987 \vol 18 \pages 545\>\~\>549
\endref

\ref\Key A2
\by \Aomoto/
\paper On the complex Selberg integral
\jour Q\.J\.Math\.Oxford \yr 1987 \vol 38 \pages 385\>\~\>399
\endref

\ref\Key AK
\by \Aomoto/ and Y\]\&Kato
\paper Gauss decomposition of connection matrices for \sym/ \$A$-type Jackson
integrals \jour \SMNS/ \vol 1 \yr 1995 \issue 4 \pages 623\>\~\>666
\endref

\ref\Key As
\by R\&Askey
\paper Some basic \hgeom/ extensions of integrals of Selberg and Andrews
\jour SIAM J.\ Math.\ Anal.{} \vol 11 \yr 1980 \pages 938\>\~\>951
\endref

\ref\Key D
\by V\]\&Dotsenko
\paper Solving the $SU(2)$ conformal field theory with the Wakimoto free-field
\rep/ \paperinfo (in Russian) \yr 1990 \pages 1\~\>31 
\endref

\ref\Key DF1
\by V\]\&Dotsenko and V\]\&Fateev
\paper Conformal algebra and multi-point correlation \fn/s in 2-D statistical
models \jour \Nucl/ \yr 1984 \vol 240 \pages 312\>\~\>348
\endref

\ref\Key DF2
\by V\]\&Dotsenko and V\]\&Fateev
\paper Four-point correlation \fn/s and the operator algebra in 2-D conformal
invariant theories with central charge $C\le 1$
\jour \Nucl/ \yr 1985, \vol 251 \pages 691\~\)734
\endref

\ref\Key EFK
\by \Etingof/, I\&Frenkel and A\&Kirillov Jr.
\book Lectures on \rep/ theory and \KZv/ \eq/s
\bookinfo Mathematical Surveys and Monographs,
\ifUS\else\ifMag\adjust{\nl}\fi\fi
\volume 58 \publ \AMSa/ \yr 1998
\endref

\ref\Key EV
\by \Etingof/ and \Varch/
\paper Dynamical Weyl groups and applications
\jour Adv\.Math. \vol 167 \yr 2002 \issue 1 \pages 74\>\~127
\endref

\ref\Key FMTV
\by \Feld/, Ya\&Markov, \VT/ and \Varch/
\paper Differential \eq/s compatible with \KZe/s
\jour Math\. Phys., Analysis and Geometry \vol 3 \yr 2000 \pages 139\>\~177
\endref

\ref\Key FSV
\by \Feld/, L\&Stevens and \Varch/
\paper Elliptic Selberg integrals and conformal blocks
\jour Preprint \yr 2002 \pages 1\~13 \info\tt math.QA/0210040
\endref

\ref\Key FV1
\by \Feld/ and \Varch/
\paper Integral \rep/s for \sol/s of the elliptic \KZvB/ \eq/
\jour \IMRN/ \yr 1995 \issue 5 \pages 221\~\>233
\endref

\ref\Key KZ
\by V\]\&Knizhnik and A\&Zamolodchikov
\paper Current algebra and Wess\)-Zumino model in two dimensions
\jour \Nucl/ \vol 247 \yr 1984 \pages 83\)\~103
\endref

\ref\Key Ma
\by I\&G\&Macdonald
\paper Some conjectures for root systems
\jour \SIAM/ \vol 13 \yr 1982 \pages 988\>\~1007
\endref

\ref\Key M
\by M\&L\&Mehta
\book Random Matrices \publ Academic Press \yr 1991
\endref

\ref\Key MaV
\by Ya\&Markov and \Varch/
\paper Solutions of \tri/ \KZe/s satisfy dynamical \deq/s
\jour Adv\. Math. \vol 166 \yr 2002 \issue 1 \pages 100\>\~147
\endref

\ref\Key MV
\by E\&Mukhin and \Varch/
\paper Remarks on critical points of phase \fn/s and norms of \Bv/s
\inbook
Adv\. Studies in Pure Math. \vol 27 \yr 2000 \pages 239\>\~\>246
\endref

\ref\Key S
\by A\&Selberg
\paper Bemerkninger om et multiplet integral
\jour Norsk Mat\.Tidsskr. \yr 1944 \vol 26 \pages 71\~\>78
\endref

\ref\Key SV
\by \SchV/ and \Varn/
\paper Arrangements of hyperplanes and Lie algebras homology
\jour \Inv/ \vol 106 \yr 1991 \pages 139\>\~194
\endref

\ref\Key TV1
\by \VT/ and \Varch/
\paper Geometry of $q$-\hgeom/ \fn/s as a bridge between Yangians and \qaff/s
\jour \Inv/ \yr 1997 \vol 128 \issue 3 \pages 501--588
\endref

\ref\Key TV2
\by \VT/ and \Varch/
\paper Difference \eq/s compatible with \tri/ \KZ/ \difl/ \eq/s
\jour \IMRN/ \yr 2000 \issue 15 \pages 801\~\)829.
\endref

\ref\Key TV3
\by \VT/ and \Varch/
\paper Duality for \KZv/ and dynamical \eq/s
\jour Acta Appl\.Math. \vol 73 \yr 2002 \issue 1\~-2 \pages 141\~154
\endref

\ref\Key V
\by \Varch/
\book Multidimensional \hgeom/ \fn/s and \rep/ theory of Lie algebras and
\qg/s \bookinfo Advanced Series in Math.\ Phys., vol\&21
\yr 1995 \publ \WSa/
\endref

\endRefs

\bye